\documentclass[12pt,twoside]{amsart}
\usepackage[latin1]{inputenc}
\usepackage{amsmath, amsthm, amscd, amsfonts, amssymb, graphicx}
\usepackage[bookmarksnumbered, plainpages]{hyperref}

\newtheorem{thm}{Theorem}[section]

\newtheorem{defn}[thm]{Definition}

\numberwithin{equation}{section}

\begin{document}

\title{\bf Geometric structures in $Sol_{3}$}
\author{Siyao Liu$^{*}$}

\thanks{{\scriptsize
\hskip -0.4 true cm \textit{2010 Mathematics Subject Classification:}
53B20; 53B25; 53B30; 53C15; 53C25; 53C35.
\newline \textit{Key words and phrases:} Semi-symmetric non-metric connection; Levi-Civita connection; $Sol_{3}$ space; Weyl conformal curvature tensor; Quasi-Einstein manifold.
\newline \textit{$^{*}$Corresponding author}}}

\maketitle

\begin{abstract}
 \indent In this paper, we conclude the geometric structures in $Sol_{3}$ with the Levi-Civita connection and the semi-symmetric non-metric connection.
 We discuss two special cases to study the geometric structures in $Sol_{3}$ with the semi-symmetric non-metric connection, namely $\widehat{P}=a(X_{3})\partial_{3}$ and $\widehat{P}=b(X_{3})\partial_{1},$ also get relatively complete conclusions respectively.
 As a special case, we investigate the geometric structures in $Sol_{3}$ with the Levi-Civita connection.
\end{abstract}

\vskip 0.2 true cm

%------------------------------------------------------------------------------------%

\pagestyle{myheadings}
\markboth{\rightline {\scriptsize  Liu}}
         {\leftline{\scriptsize }}

\bigskip
\bigskip

%------------------------------------------------------------------------------------%
%------------------------------------------------------------------------------------%

\section{ Introduction}
Based on the importance of the Morris-Thorne wormhole metric in astrophysics, Sabina in \cite{SDM} gave the geometric properties admitted by Morris-Thorne spacetime, such as Ricci generalized pseudosymmetry, Ricci generalized projectively pseudosymmetry, pseudosymmetry due to Weyl conformal curvature, semisymmetry due to conharmonic curvature, etc.
Researchers who want to investigate the geometric properties of some metrics need to describe the geometric structures of corresponding spacetimes, such as semisymmetry \cite{S1,S2,S3}, Deszcz pseudosymmetry \cite{AD}, Chaki pseudosymmetry \cite{C} and recurrent manifold \cite{R1,R2,R3}.

Agashe and Chafle introduced the notion of a semi-symmetric non-metric connection and studied some of its properties and submanifolds of a Riemannian manifold with a semi-symmetric non-metric connection \cite{AC1,AC2}.
Wang also researched non-integrable distributions with the semi-symmetric non-metric connection \cite{Wang1}.
In \cite{Wang2}, the author considered the generalized Kasner spacetimes with a semi-symmetric non-metric connection.
Liu and Wang in \cite{LW} extended the Levi-Civita connection in the definition of spacetime to the semi-symmetric non-metric connection and concluded geometric structures admitted by the Morris-Thorne wormhole metric with the semi-symmetric non-metric connection.

The space $Sol_{3}$ is a simply connected homogeneous $3$-dimensional manifold whose isometry group has dimension $3$ and it is one of the eight models of the geometry of
Thurston \cite{T}.
This space can be viewed as $\mathbb{R}^{3}$ with the metric
\begin{align}
ds^{2}=e^{2X_{3}}dX_{1}^{2}+e^{-2X_{3}}dX_{2}^{2}\pm dX_{3}^{2},
\end{align}
where $(X_{1}, X_{2}, X_{3})$ are canonical coordinates of $\mathbb{R}^{3}.$

This paper is devoted to the study of the geometric structures in $Sol_{3}$ with the semi-symmetric non-metric connection.
We select two special cases to study, namely $\widehat{P}=a(X_{3})\partial_{3}$ and $\widehat{P}=b(X_{3})\partial_{1},$ also get relatively complete conclusions respectively.
When $\widehat{P}=a(X_{3})\partial_{3},$ the space $Sol_{3}$ is projectively pseudosymmetry, Ricci generalized projectively pseudosymmetry, projectively semisymmetric, Ricci generalized projectively semisymmetric, etc.
Concurrently, it is an Einstein manifold of level $2.$
When $\widehat{P}=b(X_{3})\partial_{1},$ the space $Sol_{3}$ is a semisymmetric type manifold due to Riemann curvature, and it is an Einstein manifold of level $3.$
We investigate the geometric structures in $Sol_{3}$ with the Levi-Civita connection, that is $a(X_{3})=0$ or $b(X_{3})=0.$

A brief description of the organization of this paper is as follows.
We explain in Sec. $2$ the basic notions of several geometric structures which need to be used in this paper.
Section $3$ provides a detailed exposition of the computational process of various tensors and gives the related conclusions when $\widehat{P}=a(X_{3})\partial_{3}.$
We simultaneously treat the case when $\widehat{P}=b(X_{3})\partial_{1}.$
In Section $4,$ we detailedly describe the computational process with the Levi-Civita connection and give the geometric structures in $Sol_{3}$ with the Levi-Civita connection.
\vskip 1 true cm

\section{Preliminaries}
Let $M$ be a connected smooth semi-Riemannian manifold of dimension $n(\geq 3)$ equipped with the semi-Riemannian metric $g.$
From now on, $\nabla$ denotes the Levi-Civita connection and $X, Y, Y_{1}, Y_{2}..., U, U_{1}, U_{2}..., V, V_{1}, V_{2}..., \widehat{P}\in \mathfrak{X}(M).$
We consider the semi-symmetric non-metric connection on $M$ \cite{Wang1,Wang2}
\begin{align}
\widehat{\nabla}_{X}Y=\nabla_{X}Y+g(\widehat{P},Y)X.
\end{align}

$A$ and $E$ are two symmetric $(0,2)$-tensors, we define their Kulkarni-Nomizu product $A\wedge E$ as \cite{G1,G2,SRK}
\begin{align}
(A\wedge E)(Y_{1}, Y_{2}, U_{1}, U_{2})&=A(Y_{1}, U_{2})E(Y_{2},U_{1})-A(Y_{1},U_{1})E(Y_{2}, U_{2})\\
&+A(Y_{2},U_{1})E(Y_{1}, U_{2})-A(Y_{2}, U_{2})E(Y_{1}, U_{1}).\nonumber
\end{align}
For a symmetric $(0,2)$-tensor $A,$ we can define the endomorphism $Y_{1}\wedge_{A} Y_{2}$ as \cite{DDHKS,SH1,SH2}
\begin{align}
(Y_{1}\wedge_{A} Y_{2})U=A(Y_{2}, U)Y_{1}-A(Y_{1},U)Y_{2}.
\end{align}

Write
\begin{align}
\mathcal{E}_{R}&(Y_{1}, Y_{2})=\widehat{\nabla}_{Y_{1}}\widehat{\nabla}_{Y_{2}}-\widehat{\nabla}_{Y_{2}}\widehat{\nabla}_{Y_{1}}-\widehat{\nabla}_{[Y_{1}, Y_{2}]}.
\end{align}
We introduce four endomorphisms on $M$ \cite{K,SC}
\begin{align}
\mathcal{E}_{C}&=\mathcal{E}_{R}-\frac{1}{n-2}(\mathcal{J}\wedge_{g}+\wedge_{g}\mathcal{J}-\frac{\kappa}{n-1}\wedge_{g}),\\
\mathcal{E}_{K}&=\mathcal{E}_{R}-\frac{1}{n-2}(\mathcal{J}\wedge_{g}+\wedge_{g}\mathcal{J}),\\
\mathcal{E}_{W}&=\mathcal{E}_{R}-\frac{\kappa}{n(n-1)}\wedge_{g},\\
\mathcal{E}_{P}&=\mathcal{E}_{R}-\frac{1}{n-1}\wedge_{Ric},
\end{align}
where $\mathcal{J}$ is the Ricci operator defined as $Ric(Y_{1}, Y_{2})=g(Y_{1}, \mathcal{J}Y_{2}),$
the Ricci tensor $Ric$ is defined by $Ric(Y_{1}, Y_{2})=tr\{X\rightarrow \mathcal{E}_{R}(X, Y_{1})Y_{2}\}$
and $\kappa$ denotes the scalar curvature with the semi-symmetric non-metric connection.
For an endomorphism $\mathcal{E}(U_{1}, U_{2})$ one can define the following $(0,4)$ type tensor field as
\begin{align}
E(U_{1}, U_{2}, U_{3}, U_{4})=g(\mathcal{E}(U_{1}, U_{2})U_{3}, U_{4}).
\end{align}
If we replace $\mathcal{E}$ with the endomorphisms $\mathcal{E}_{R}$ $\mathcal{E}_{C}$ $\mathcal{E}_{K}$ $\mathcal{E}_{W}$ and $\mathcal{E}_{P}$ above, then we can  obtain the $(0,4)$ type Riemann curvature $R,$ the Weyl conformal curvature $C,$ the conharmonic curvature $K,$ the concircular curvature $W$ and the projective curvature $P$ about $\widehat{\nabla}.$
These tensors locally are given by
\begin{align}
R_{hkij}&=g_{h\alpha}\left(\partial_{i}\widehat{\Gamma}^{\alpha}_{kj}-\partial_{j}\widehat{\Gamma}^{\alpha}_{ki}+\widehat{\Gamma}^{\beta}_{kj}\widehat{\Gamma}^{\alpha}_{\beta i}-\widehat{\Gamma}^{\beta}_{ki}\widehat{\Gamma}^{\alpha}_{\beta j}\right),\\
C_{hkij}&=\left(R-\frac{1}{n-2}(g\wedge Ric)+\frac{\kappa}{2(n-1)(n-2)}(g\wedge g)\right)_{hkij},\\
K_{hkij}&=\left(R-\frac{1}{n-2}(g\wedge Ric)\right)_{hkij},\\
W_{hkij}&=\left(R-\frac{\kappa}{2n(n-1)}(g\wedge g)\right)_{hkij},\\
P_{hkij}&=R_{hkij}-\frac{1}{n-1}(g_{hi}Ric_{kj}-g_{ki}Ric_{hj}),
\end{align}
where $g_{h\alpha}$ is the metric matrix, $\partial_{i}$ is a natural local frame, $\widehat{\Gamma}^{\alpha}_{kj}$ is the Christofel coefficient of the semi-symmetric non-metric connection.

Here $F$ is a $(0,k),$ $k\geq 1,$ type tensor field, we have a $(0,k+2)$ type tensor field $E\cdot F$ as follows \cite{DG,SH3}:
\begin{align}
(E\cdot F)(Y_{1}, Y_{2},..., Y_{k}, U_{1}, U_{2})&=(\mathcal{E}(U_{1}, U_{2})\cdot F)(Y_{1}, Y_{2},..., Y_{k})\\
&=-F(\mathcal{E}(U_{1}, U_{2})Y_{1}, Y_{2},..., Y_{k})-...-F(Y_{1}, Y_{2},..., \mathcal{E}(U_{1}, U_{2})Y_{k}),\nonumber
\end{align}
and a $(0,k+2)$ type tensor $Q(E,F),$ named Tachibana tensor, as follows \cite{T2,DGPSS}:
\begin{align}
&Q(E, F)(Y_{1}, Y_{2},..., Y_{k}, U_{1}, U_{2})=((U_{1}\wedge_{E} U_{2})\cdot F)(Y_{1}, Y_{2},..., Y_{k})\\
&=E(U_{1}, Y_{1})F(U_{2}, Y_{2},..., Y_{k})+...+E(U_{1}, Y_{k})F(Y_{1}, Y_{2},..., U_{2})\nonumber\\
&-E(U_{2}, Y_{1})F(U_{1}, Y_{2},..., Y_{k})-...-E(U_{2}, Y_{k})F(Y_{1}, Y_{2},..., U_{1}).\nonumber
\end{align}
\begin{defn} \cite{R4,R5,SH4,SH5}
If the tensor $E\cdot F$ is linearly dependent with $Q(Z,F)$ then it is called as $F$-pseudosymmetric type manifold due to $E,$  if $E\cdot F=0$ then it is called $F$-semisymmetric type manifold due to $E.$
\end{defn}
In particular, if $E\cdot F= f_{F}Q(Ric,F)$ ($f_{F}$ being some smooth scalar function on $M$), then it is called a Ricci generalized $F$-pseudosymmetric type manifold due to $E.$
A $F$-pseudosymmetric manifold due to $E$ is said to be a Deszcz pseudosymmetric manifold if $E=R,$ $F=R$ and $Z=g,$ and a $F$-semisymmetric type manifold due to $E$ is said to be a semisymmetric manifold if $E=R,$ $F=R.$
Again for $E=C,$ $F=C$ and $Z=g$ the manifold is said to have a pseudosymmetric Weyl conformal curvature tensor.

With the definition in \cite{S,SYH}, we have
\begin{defn}
For a scalar $\alpha$ and $0\leq k\leq n-1$, if $rank(Ric-\alpha g)=k,$ then $M$ is called as a $k$-quasi-Einstein manifold and for $k=1$ (resp., $k=0$) it is a quasi-Einstein (resp., Einstein) manifold. A Ricci simple manifold is a quasi-Einstein manifold for $\alpha = 0$.
\end{defn}
Our definition is similar to the one given in \cite{B,D,DGJZ}.
\begin{defn}
A generalized Roter type manifold $M$ is a manifold whose Riemann curvature can explicitly be expressed as
\begin{align}
R=(\mu_{11}Ric^{2}+\mu_{12}Ric+\mu_{13}g)\wedge Ric^{2}+(\mu_{22}Ric+\mu_{23}g)\wedge Ric+\mu_{33}(g\wedge g),
\end{align}
where $\mu_{ij}$ are some scalars, the Ricci tensor of level $k$ is defined by $Ric^{k}(Y_{1}, Y_{2})=g(Y_{1}, \mathcal{J}^{k-1}Y_{2}).$ If such linear dependency reduces to the linear dependency of $R,$ $g\wedge g,$ $g\wedge Ric,$ $Ric\wedge Ric$ then it is called a manifold of Roter.
\end{defn}
Since the definition in \cite{SH6}, Einstein manifold of level $2$ (resp., $3$ and $4$) about the connection $\widehat{\nabla}$ is well defined.
\begin{defn} \cite{SH6}
A semi-Riemannian manifold $M$ is said to be Einstein manifold of level $2$ (resp., $3$ and $4$) about the connection $\widehat{\nabla}$ if
\begin{align}
&Ric^{2}+\lambda_{1}Ric+\lambda_{2}g=0,\\
&Ric^{3}+\lambda_{3}Ric^{2}+\lambda_{4}Ric+\lambda_{5}g=0,\\
&Ric^{4}+\lambda_{6}Ric^{3}+\lambda_{7}Ric^{2}+\lambda_{8}Ric+\lambda_{9}g=0,
\end{align}
where $\lambda_{i}(1\leq\lambda_{i}\leq9)$ are some scalars.
\end{defn}
Based on the concept of \cite{G,SJ}, we have
\begin{defn}
The Ricci tensor of the semi-Riemannian manifold $M$ is said to be cyclic parallel if the relation
\begin{align}
\underset{Y_{1}, Y_{2}, Y_{3}}{\mathcal{S}} (\widehat{\nabla}_{Y_{1}}Ric)(Y_{2}, Y_{3})=0
\end{align}
holds and it is said to be Codazzi if $(\widehat{\nabla}_{Y_{1}}Ric)(Y_{2}, Y_{3})=(\widehat{\nabla}_{Y_{2}}Ric)(Y_{1}, Y_{3})$ holds where $\mathcal{S}$ is a cyclic sum over $Y_{1},$ $Y_{2}$ and
$Y_{3}.$
\end{defn}
\begin{defn} \cite{MM1,MM2,MM3}
A symmetric $(0,2)$-tensor $Z$ on a semi-Riemannian M is said to be $T$-compatible if
\begin{align}
\underset{Y_{1}, Y_{2}, Y_{3}}{\mathcal{S}} T(\mathcal{Z}Y_{1}, U, Y_{2}, Y_{3})=0
\end{align}
holds, where $\mathcal{Z}$ is the endomorphism corresponding to $Z$ defined as $g(\mathcal{Z}Y_{1}, Y_{2})= Z(Y_{1}, Y_{2}).$
\end{defn}
Replacing $T$ by the curvatures $R,$ $C,$ $W,$ $P$ and $K$ we can define the corresponding curvature compatibilities.

We develop the definition of \cite{MS1,MS2}, then  we could find
\begin{defn}
The curvature $2$-forms $\Omega^{M}_{(E)l}$ are recurrent if and only if
\begin{align}
\underset{Y_{1}, Y_{2}, Y_{3}}{\mathcal{S}} (\widehat{\nabla}_{Y_{1}}E)(Y_{2}, Y_{3}, U, Y)=\underset{Y_{1}, Y_{2}, Y_{3}}{\mathcal{S}} \sigma(Y_{1})E(Y_{2}, Y_{3}, U, Y).
\end{align}
\end{defn}
\vskip 1 true cm

\section{The geometric structures in $Sol_{3}$ with the semi-symmetric non-metric connection}
The components of the metric $(1.1)$ are given by
\begin{align}
g_{11}=e^{2 X_3};~g_{22}=e^{-2 X_3};~g_{33}=\pm 1;~g_{ij}=0,otherwise.
\end{align}

Let $\widehat{P}=a(X_{3})\frac{\partial}{\partial X_{3}}=a(X_{3})\partial_{3},$ where $a(X_{3})$ is a continuous function.
At this time, the non-vanishing components of the Christofel coefficient $\widehat{\Gamma}^{\alpha}_{kj}$ are given below
\begin{align}
&\widehat{\Gamma}^{1}_{13}=1 \pm a(X_{3}),~\widehat{\Gamma}^{1}_{31}=1,~\widehat{\Gamma}^{2}_{23}=-1 \pm a(X_{3}),~\widehat{\Gamma}^{2}_{32}=-1,\\
&\widehat{\Gamma}^{3}_{11}=\mp e^{2 X_3},~\widehat{\Gamma}^{3}_{22}=\pm e^{-2 X_3},~\widehat{\Gamma}^{3}_{33}=\pm a(X_{3}).\nonumber
\end{align}
Then the non-zero components of the Riemann curvature tensor $R,$ Ricci tensor $Ric$ and scalar curvature $\kappa$ of (1.1) are given by
\begin{align}
&R_{1212}=\pm 1,~R_{1313}=-e^{2 X_3},~R_{2323}=-e^{-2 X_3};\\
&Ric_{33}=2;~\kappa=\pm 2.
\end{align}
Thus, the non-zero components of $\widehat{\nabla} R$ and $\widehat{\nabla} Ric$ are calculated as below
\begin{align}
&\widehat{\nabla}_{2} R_{1213}=\pm 2-a(X_{3}),~\widehat{\nabla}_{1} R_{1223}=\pm 2+a(X_{3}),\\
&\widehat{\nabla}_{3} R_{1313}=\pm 2e^{2 X_3}a(X_{3}),~\widehat{\nabla}_{3} R_{2323}=\pm 2e^{-2 X_3}a(X_{3});\nonumber\\
&\widehat{\nabla}_{1} Ric_{13}=\pm 2e^{2 X_3},~\widehat{\nabla}_{2} Ric_{23}=\mp 2e^{-2 X_3},~\widehat{\nabla}_{3} Ric_{33}=\mp 4a(X_{3}).
\end{align}
According to some computations about $g\wedge Ric$ and $g\wedge g$, we can get the non zero components of the tensor $C,$ $K,$ $W,$ $P$ as follows
\begin{align}
&K_{1212}=\pm 1,~K_{1313}=e^{2 X_3},~K_{2323}=e^{-2 X_3};\\
&W_{1212}=\pm \frac{4}{3},~W_{1313}=-\frac{2}{3}e^{2 X_3},~W_{2323}=-\frac{2}{3}e^{-2 X_3};\\
&P_{1212}=-P_{1221}=\pm 1,~P_{1313}=-e^{2 X_3},~P_{2323}=-e^{-2 X_3}.
\end{align}
It is easy to check that
\begin{align}
&\widehat{\nabla}_{2} K_{1213}=\widehat{\nabla}_{2} K_{1312}=-a(X_{3}),~\widehat{\nabla}_{1} K_{1223}=\widehat{\nabla}_{1} K_{2312}=a(X_{3}),\\
&\widehat{\nabla}_{3} K_{1313}=\mp 2e^{2 X_3}a(X_{3}),~\widehat{\nabla}_{3} K_{2323}=\mp 2e^{-2 X_3}a(X_{3});\nonumber\\
&\widehat{\nabla}_{2} W_{1213}=\widehat{\nabla}_{2} W_{1312}=\pm 2-\frac{4}{3}a(X_{3}),~\widehat{\nabla}_{1} W_{1223}=\widehat{\nabla}_{1} W_{2312}=\pm 2+\frac{4}{3}a(X_{3}),\\
&\widehat{\nabla}_{3} W_{1313}=\pm \frac{4}{3}e^{2 X_3}a(X_{3}),~\widehat{\nabla}_{3} W_{2323}=\pm \frac{4}{3}e^{-2 X_3}a(X_{3});\nonumber\\
&\widehat{\nabla}_{2} P_{1213}=\widehat{\nabla}_{2} P_{1312}=\pm 2-a(X_{3}),~\widehat{\nabla}_{1} P_{1223}=-\widehat{\nabla}_{1} P_{2321}=\pm 2+a(X_{3}),\\
&\widehat{\nabla}_{2} P_{1231}=\widehat{\nabla}_{2} P_{1321}=\mp 1+a(X_{3}),~\widehat{\nabla}_{1} P_{1232}=-\widehat{\nabla}_{1} P_{2312}=\mp 1-a(X_{3}),\nonumber\\
&\widehat{\nabla}_{1} P_{1311}=\mp e^{4 X_3},~\widehat{\nabla}_{3} P_{1313}=\pm 2e^{2 X_3}a(X_{3}),~\widehat{\nabla}_{1} P_{1333}=e^{2 X_3}(1\pm a(X_{3})),~\nonumber\\
&\widehat{\nabla}_{2} P_{2322}=\pm e^{-4 X_3},~\widehat{\nabla}_{3} P_{2323}=\pm 2e^{-2 X_3}a(X_{3}),~\widehat{\nabla}_{2} P_{2333}=e^{-2 X_3}(-1\pm  a(X_{3})).\nonumber
\end{align}
Again the non-zero components of the tensor $R\cdot R,$ $R\cdot C,$ $R\cdot K,$ $R\cdot W,$ $R\cdot P$ are written below
\begin{align}
&(R\cdot R)_{121323}=\mp 2,~(R\cdot R)_{122313}=(R\cdot R)_{132312}=\pm 2;\\
&(R\cdot W)_{121323}=(R\cdot W)_{131223}=\mp 2,~(R\cdot W)_{122313}=(R\cdot W)_{231213}=\pm 2;\\
&(R\cdot P)_{121323}=(R\cdot P)_{131223}=(R\cdot P)_{232113}=\mp 2,~(R\cdot P)_{122313}=\pm 2,\\
&(R\cdot P)_{123123}=(R\cdot P)_{132123}=(R\cdot P)_{231213}=\pm 1,~(R\cdot P)_{123213}=\mp 1,\nonumber\\
&(R\cdot P)_{131113}=\mp e^{4 X_3},~(R\cdot P)_{133313}=e^{2 X_3}.\nonumber
\end{align}
An easy computation shows that
\begin{align}
&(K\cdot R)_{121323}=(K\cdot R)_{131223}=\pm 2,~(K\cdot R)_{122313}=(K\cdot R)_{231213}=\mp 2;\\
&(K\cdot W)_{121323}=(K\cdot W)_{131223}=\pm 2,~(K\cdot W)_{122313}=(K\cdot W)_{231213}=\mp 2;\\
&(K\cdot P)_{121323}=(K\cdot P)_{131223}=(K\cdot P)_{232113}=\pm 2,~(K\cdot P)_{122313}=\mp 2,\\
&(K\cdot P)_{123123}=(K\cdot P)_{132123}=(K\cdot P)_{231213}=\mp 1,~(K\cdot P)_{123213}=\pm 1,\nonumber\\
&(K\cdot P)_{131113}=\pm e^{4 X_3},~(K\cdot P)_{133313}=-e^{2 X_3}.\nonumber
\end{align}
In addition, we have
\begin{align}
&(W\cdot R)_{121323}=(W\cdot R)_{131223}=\mp \frac{4}{3},~(W\cdot R)_{122313}=(W\cdot R)_{231213}=\pm \frac{4}{3};\\
&(W\cdot W)_{121323}=(W\cdot W)_{131223}=\mp \frac{4}{3},~(W\cdot W)_{122313}=(W\cdot W)_{231213}=\pm \frac{4}{3};\\
&(W\cdot P)_{121323}=(W\cdot P)_{131223}=(W\cdot P)_{232113}=\mp \frac{4}{3},~(W\cdot P)_{122313}=\pm \frac{4}{3},\\
&(W\cdot P)_{123123}=(W\cdot P)_{132123}=(W\cdot P)_{231213}=\pm \frac{2}{3},~(W\cdot P)_{123213}=\mp \frac{2}{3},\nonumber\\
&(W\cdot P)_{131113}=\mp \frac{2}{3}e^{4 X_3},~(W\cdot P)_{133313}=\frac{2}{3}e^{2 X_3}\nonumber
\end{align}
and
\begin{align}
&(P\cdot R)_{121323}=(P\cdot R)_{131223}=\mp 1,~(P\cdot R)_{122313}=(P\cdot R)_{231213}=\pm 1;\\
&(P\cdot K)_{121323}=(P\cdot K)_{131223}=\pm 1,~(P\cdot K)_{122313}=(P\cdot K)_{231213}=\mp 1;\\
&(P\cdot W)_{121323}=(P\cdot W)_{131223}=\mp \frac{2}{3},~(P\cdot W)_{122313}=(P\cdot W)_{231213}=\pm \frac{2}{3};\\
&(P\cdot P)_{121323}=(P\cdot P)_{131223}=(P\cdot P)_{232113}=\mp 1,~(P\cdot P)_{122313}=\pm 1,\\
&(P\cdot P)_{131113}=\mp e^{4 X_3}.\nonumber
\end{align}
By $(2.16),$ the non-vanishing components of the Tachibana tensor are given by
\begin{align}
&Q(Ric, R)_{121323}=Q(Ric, R)_{131223}=\mp 2,~Q(Ric, R)_{122313}=Q(Ric, R)_{231213}=\pm 2;\\
&Q(Ric, K)_{121323}=Q(Ric, K)_{131223}=\mp 2,~Q(Ric, K)_{122313}=Q(Ric, K)_{231213}=\pm 2;\\
&Q(Ric, W)_{121323}=Q(Ric, W)_{131223}=\mp \frac{8}{3},~Q(Ric, W)_{122313}=Q(Ric, W)_{231213}=\pm \frac{8}{3};\\
&Q(Ric, P)_{121323}=Q(Ric, P)_{131223}=\mp 2,~Q(Ric, P)_{122313}=Q(Ric, P)_{231213}=\pm 2,\\
&Q(Ric, P)_{133313}=2e^{2 X_3}.\nonumber
\end{align}
We see at once that
\begin{align}
&Q(g, R)_{121323}=Q(g, R)_{131223}=-2,~Q(g, R)_{122313}=Q(g, R)_{231213}=2;\\
&Q(g, W)_{121323}=Q(g, W)_{131223}=-2,~Q(g, W)_{122313}=Q(g, W)_{231213}=2;\\
&Q(g, P)_{121323}=Q(g, P)_{131223}=Q(g, P)_{232113}=-2,~Q(g, P)_{122313}=2,\\
&Q(g, P)_{123123}=Q(g, P)_{132123}=Q(g, P)_{231213}=1,~Q(g, P)_{123213}=-1,\nonumber\\
&Q(g, P)_{131113}=-e^{4 X_3},~Q(g, P)_{133313}=\pm e^{2 X_3}.\nonumber
\end{align}

We can now formulate our first results.
\begin{thm}
The metric $(1.1)$ fulfills the following geometric structures about the connection $\widehat{\nabla}$:

(I)(1)$C\cdot R=0$ from this it is the Riemann curvature semisymmetric type manifold due to the Weyl conformal curvature;\\
(2)$K\cdot R=-Q(Ric, R)$ hence it is Ricci generalized the Riemann curvature pseudosymmetric type manifold due to the conharmonic curvature;\\
(3)$W\cdot R=\frac{2}{3}Q(Ric, R)$ thus it is Ricci generalized the Riemann curvature pseudosymmetric type manifold due to the concircular curvature;\\
(4)Ricci generalized projectively pseudosymmetric because $P\cdot R=\frac{1}{2}Q(Ric, R)$;\\
(5)$K\cdot R=\mp Q(g, R)$ so it is the Riemann curvature pseudosymmetric type manifold due to the conharmonic curvature;\\
(6)$W\cdot R=\pm \frac{2}{3}Q(g, R)$ from this it is the Riemann curvature pseudosymmetric type manifold due to the concircular curvature;\\
(7)$P\cdot R=\pm \frac{1}{2}Q(g, R)$ i.e. it is the Riemann curvature pseudosymmetric type manifold due to the projective curvature;\\
(8)$R\cdot C=0$ then it is the Weyl conformal curvature semisymmetric type manifold due to the Riemann curvature;\\
(9)$C\cdot C=0$ for this reason it is the Weyl conformal curvature semisymmetric type manifold due to the Weyl conformal curvature;\\
(10)$K\cdot C=0$ in this way it is the Weyl conformal curvature semisymmetric type manifold due to the conharmonic curvature;\\
(11)$W\cdot C=0$ hence it is the Weyl conformal curvature semisymmetric type manifold due to the concircular curvature;\\
(12)$P\cdot C=0$ so that it is the Weyl conformal curvature semisymmetric type manifold due to the projective curvature;\\
(13)$R\cdot K=0$ i.e. it is the conharmonic curvature semisymmetric type manifold due to the Riemann curvature;\\
(14)$C\cdot K=0$ for this reason it is the conharmonic curvature semisymmetric type manifold due to the Weyl conformal curvature;\\
(15)$K\cdot K=0$ so that it is the conharmonic curvature semisymmetric type manifold due to the conharmonic curvature;\\
(16)$W\cdot K=0$ in this way it is the conharmonic curvature semisymmetric type manifold due to the concircular curvature;\\
(17)$P\cdot K=-\frac{1}{2}Q(Ric, K)$ for this reason it is Ricci generalized the conharmonic curvature pseudosymmetric type manifold due to the projective curvature;\\
(18)$R\cdot W=\frac{3}{4}Q(Ric, W)$ hence it is Ricci generalized the concircular curvature pseudosymmetric type manifold due to the Riemann curvature;\\
(19)$K\cdot W=-\frac{3}{4}Q(Ric, W)$ from this it is Ricci generalized the concircular curvature pseudosymmetric type manifold due to the conharmonic curvature;\\
(20)$W\cdot W=\frac{1}{2}Q(Ric, W)$ therefore it is Ricci generalized the concircular curvature pseudosymmetric type manifold due to the concircular curvature;\\
(21)$P\cdot W=\frac{1}{4}Q(Ric, W)$ so it is Ricci generalized the concircular curvature pseudosymmetric type manifold due to the projective curvature;\\
(22)$R\cdot W=\pm Q(g, W)$ then it is the concircular curvature pseudosymmetric type manifold due to the Riemann curvature;\\
(23)Concircular curvature pseudosymmetric type manifold due to the conharmonic curvature i.e. $K\cdot W=\mp Q(g, W)$;\\
(24)$W\cdot W=\pm \frac{2}{3}Q(g, W)$ in this way it is the concircular curvature pseudosymmetric type manifold due to the concircular curvature;\\
(25)$P\cdot W=\pm \frac{1}{3}Q(g, W)$ thus it is the concircular curvature pseudosymmetric type manifold due to the projective curvature;\\
(26)$R\cdot P=\pm Q(g, P)$ therefore it is the projective curvature pseudosymmetric type manifold due to the Riemann curvature;\\
(27)$K\cdot P=\mp Q(g, P)$ for this reason it is the projective curvature pseudosymmetric type manifold due to the conharmonic curvature;\\
(28)$W\cdot P=\pm \frac{2}{3}Q(g, P)$  hence it is the projective curvature pseudosymmetric type manifold due to the concircular curvature.

(II)$Ric=\alpha(\eta\bigotimes\eta)$ for $\alpha=2$ and $\eta=\{0, 0, 1\},$ thus it is a Ricci simple manifold and a $2$-quasi-Einstein manifold as $rank(Ric-\alpha g)=2$ when $g_{33}=1,$ or a $3$-quasi-Einstein manifold as $rank(Ric-\alpha g)=3$ when $g_{33}=-1.$

(III)The manifold is the generalized Roter type manifold since $R=g\wedge Ric\mp \frac{1}{2} g\wedge g.$

(IV)The condition $Ric^{2}\mp 2 Ric=0$ holds, hence it is a Einstein manifold of level $2.$

(V)Ricci tensor is neither of Codazzi type nor cyclic parallel.

(VI)Ricci tensor is Riemann compatible, Weyl conformal compatible, conharmonic compatible and concircular compatible.

(VII)(1)The Riemann curvature are recurrent for the $1$-form $\{0, 0, \pm2a(X_{3})\};$\\
(2)The conharmonic curvature are recurrent for the $1$-form $\{0, 0, \pm2a(X_{3})\};$\\
(3)The concircular curvature are recurrent for the $1$-form $\{0, 0, \pm2a(X_{3})\}.$
\end{thm}

Consider $\widehat{P}=b(X_{3})\frac{\partial}{\partial X_{1}}=b(X_{3})\partial_{1},$ where $b(X_{3})$ is a continuous function. We thus get
\begin{align}
&\widehat{\Gamma}^{1}_{11}=\widehat{\Gamma}^{2}_{21}=\widehat{\Gamma}^{3}_{31}=e^{2 X_3}b(X_{3}),~\widehat{\Gamma}^{3}_{11}=\mp e^{2 X_3},\\
&\widehat{\Gamma}^{1}_{13}=\widehat{\Gamma}^{1}_{31}=1,~\widehat{\Gamma}^{3}_{22}=\pm e^{-2 X_3},~\widehat{\Gamma}^{2}_{23}=\widehat{\Gamma}^{2}_{32}=-1.\nonumber
\end{align}
By some computations, we can get
\begin{align}
&R_{1113}=-e^{4 X_3}(2b(X_{3})+b'(X_{3})),~R_{1212}=\pm 1,~R_{1313}=-e^{2 X_3},\\
&R_{2213}=-(2b(X_{3})+b'(X_{3})),~R_{2323}=-e^{-2 X_3},~R_{3313}=\mp e^{2 X_3}(2b(X_{3})+b'(X_{3}));\nonumber\\
&Ric_{13}=-Ric_{31}=e^{2 X_3}(2b(X_{3})+b'(X_{3})),~Ric_{33}=2;~\kappa=\pm 2,
\end{align}
where $b'(X_{3})=\frac{\partial b(X_{3})}{\partial X_{3}},$ $b''(X_{3})=\frac{\partial b'(X_{3})}{\partial X_{3}}.$
We notice that
\begin{align}
&C_{1113}=-e^{4 X_3}(2b(X_{3})+b'(X_{3})),~C_{1223}=-(2b(X_{3})+b'(X_{3})),\\
&C_{2213}=-(2b(X_{3})+b'(X_{3})),~C_{2312}=2b(X_{3})+b'(X_{3}),~C_{3313}=\mp e^{2 X_3}(2b(X_{3})+b'(X_{3}));\nonumber\\
&K_{1113}=-e^{4 X_3}(2b(X_{3})+b'(X_{3})),~K_{1212}=\pm 1,~K_{1223}=-(2b(X_{3})+b'(X_{3})),\\
&K_{1313}=e^{2 X_3},~K_{2213}=-(2b(X_{3})+b'(X_{3})),~K_{2312}=2b(X_{3})+b'(X_{3}),\nonumber\\
&K_{2323}=e^{-2 X_3},~K_{3313}=\mp e^{2 X_3}(2b(X_{3})+b'(X_{3}));\nonumber\\
&W_{1113}=-e^{4 X_3}(2b(X_{3})+b'(X_{3})),~W_{1212}=\pm \frac{4}{3},~W_{1313}=-\frac{2}{3}e^{2 X_3},\\
&W_{2213}=-(2b(X_{3})+b'(X_{3})),~W_{2323}=- \frac{2}{3}e^{-2 X_3},~W_{3313}=\mp e^{2 X_3}(2b(X_{3})+b'(X_{3}));\nonumber\\
&P_{1113}=-P_{1131}=-e^{4 X_3}(2b(X_{3})+b'(X_{3})),~P_{1212}=-P_{1221}=\pm 1,\\
&P_{1232}=\frac{1}{2}(2b(X_{3})+b'(X_{3})),~P_{1311}=\frac{1}{2}e^{4 X_3}(2b(X_{3})+b'(X_{3})),\nonumber\\
&P_{1313}=-e^{2 X_3},~P_{1333}=\pm\frac{1}{2}e^{2 X_3}(2b(X_{3})+b'(X_{3})),\nonumber\\
&P_{2213}=-P_{2231}=-(2b(X_{3})+b'(X_{3})),~P_{2312}=\frac{1}{2}(2b(X_{3})+b'(X_{3})),\nonumber\\
&P_{2323}=-e^{-2 X_3},~P_{3313}=-P_{3331}=\mp e^{2 X_3}(2b(X_{3})+b'(X_{3})).\nonumber
\end{align}
Obviously, we can get
\begin{align}
&(R\cdot R)_{111223}=\mp e^{2 X_3}(2b(X_{3})+b'(X_{3})),~(R\cdot R)_{121323}=\mp 2,\\
&(R\cdot R)_{112312}=\mp e^{2 X_3}(2b(X_{3})+b'(X_{3})),~(R\cdot R)_{122313}=\pm 2,\nonumber\\
&(R\cdot R)_{122311}=\pm 2e^{2 X_3}(2b(X_{3})+b'(X_{3})),~(R\cdot R)_{122322}=\pm 2e^{-2 X_3}(2b(X_{3})+b'(X_{3})),\nonumber\\
&(R\cdot R)_{122333}=2(2b(X_{3})+b'(X_{3})),~(R\cdot R)_{131223}=\mp 2,\nonumber
\end{align}
\begin{align}
&(R\cdot R)_{221223}=\mp e^{-2 X_3}(2b(X_{3})+b'(X_{3})),~(R\cdot R)_{222312}=\mp e^{-2 X_3}(2b(X_{3})+b'(X_{3})),\nonumber\\
&(R\cdot R)_{231211}=\pm 2e^{2 X_3}(2b(X_{3})+b'(X_{3})),~(R\cdot R)_{231213}=\pm 2,\nonumber\\
&(R\cdot R)_{231222}=\pm 2e^{-2 X_3}(2b(X_{3})+b'(X_{3})),~(R\cdot R)_{231233}=2(2b(X_{3})+b'(X_{3})),\nonumber\\
&(R\cdot R)_{331223}=-(2b(X_{3})+b'(X_{3})),~(R\cdot R)_{332312}=-(2b(X_{3})+b'(X_{3}));\nonumber\\
&(R\cdot C)_{111223}=\mp e^{2 X_3}(2b(X_{3})+b'(X_{3})),~(R\cdot C)_{112312}=\mp e^{2 X_3}(2b(X_{3})+b'(X_{3})),\\
&(R\cdot C)_{121312}=\pm e^{2 X_3}(2b(X_{3})+b'(X_{3})),~(R\cdot C)_{131212}=\mp e^{2 X_3}(2b(X_{3})+b'(X_{3})),\nonumber\\
&(R\cdot C)_{132323}=2b(X_{3})+b'(X_{3}),~(R\cdot C)_{221223}=\mp e^{-2 X_3}(2b(X_{3})+b'(X_{3})),\nonumber\\
&(R\cdot C)_{222312}=\mp e^{-2 X_3}(2b(X_{3})+b'(X_{3})),~(R\cdot C)_{231323}=-(2b(X_{3})+b'(X_{3})),\nonumber\\
&(R\cdot C)_{331223}=-(2b(X_{3})+b'(X_{3})),~(R\cdot C)_{332312}=-(2b(X_{3})+b'(X_{3}));\nonumber\\
&(R\cdot W)_{111223}=\mp e^{2 X_3}(2b(X_{3})+b'(X_{3})),~(R\cdot W)_{121323}=\mp 2,\\
&(R\cdot W)_{112312}=\mp e^{2 X_3}(2b(X_{3})+b'(X_{3})),~(R\cdot W)_{122313}=\pm 2,\nonumber\\
&(R\cdot W)_{122311}=\pm 2e^{2 X_3}(2b(X_{3})+b'(X_{3})),~(R\cdot W)_{122322}=\pm 2e^{-2 X_3}(2b(X_{3})+b'(X_{3})),\nonumber\\
&(R\cdot W)_{122333}=2(2b(X_{3})+b'(X_{3})),~(R\cdot W)_{131223}=\mp 2,\nonumber\\
&(R\cdot W)_{221223}=\mp e^{-2 X_3}(2b(X_{3})+b'(X_{3})),~(R\cdot W)_{222312}=\mp e^{-2 X_3}(2b(X_{3})+b'(X_{3})),\nonumber\\
&(R\cdot W)_{231211}=\pm 2e^{2 X_3}(2b(X_{3})+b'(X_{3})),~(R\cdot W)_{231213}=\pm 2,\nonumber\\
&(R\cdot W)_{231222}=\pm 2e^{-2 X_3}(2b(X_{3})+b'(X_{3})),~(R\cdot W)_{231233}=2(2b(X_{3})+b'(X_{3})),\nonumber\\
&(R\cdot W)_{331223}=-(2b(X_{3})+b'(X_{3})),~(R\cdot W)_{332312}=-(2b(X_{3})+b'(X_{3}));\nonumber\\
&(R\cdot P)_{111223}=-(R\cdot P)_{112123}=\mp e^{2 X_3}(2b(X_{3})+b'(X_{3})),\\
&(R\cdot P)_{112312}=-(R\cdot P)_{113212}=\mp e^{2 X_3}(2b(X_{3})+b'(X_{3})),\nonumber\\
&(R\cdot P)_{121123}=-(R\cdot P)_{211123}=\pm \frac{1}{2}e^{2 X_3}(2b(X_{3})+b'(X_{3})),\nonumber\\
&(R\cdot P)_{121323}=-(R\cdot P)_{211323}=\mp 2,~(R\cdot P)_{122311}=\pm 2e^{2 X_3}(2b(X_{3})+b'(X_{3})),\nonumber\\
&(R\cdot P)_{122223}=-(R\cdot P)_{212223}=\pm \frac{1}{2}e^{-2 X_3}(2b(X_{3})+b'(X_{3})),\nonumber\\
&(R\cdot P)_{122313}=-(R\cdot P)_{212313}=\pm 2,~(R\cdot P)_{122322}=\pm 2e^{-2 X_3}(2b(X_{3})+b'(X_{3})),\nonumber\\
&(R\cdot P)_{122333}=-(R\cdot P)_{212333}=2(2b(X_{3})+b'(X_{3})),\nonumber\\
&(R\cdot P)_{123112}=-(R\cdot P)_{213112}=\mp \frac{1}{2}e^{2 X_3}(2b(X_{3})+b'(X_{3})),\nonumber\\
&(R\cdot P)_{123123}=-(R\cdot P)_{213123}=\pm 1,~(R\cdot P)_{123211}=\mp e^{2 X_3}(2b(X_{3})+b'(X_{3})),\nonumber\\
&(R\cdot P)_{123213}=-(R\cdot P)_{213213}=\mp 1,~(R\cdot P)_{123222}=\mp e^{-2 X_3}(2b(X_{3})+b'(X_{3})),\nonumber\\
&(R\cdot P)_{123233}=-(2b(X_{3})+b'(X_{3})),~(R\cdot P)_{131111}=\mp e^{6 X_3}(2b(X_{3})+b'(X_{3})),\nonumber
\end{align}
\begin{align}
&(R\cdot P)_{131113}=-(R\cdot P)_{311113}=\mp e^{4 X_3},~(R\cdot P)_{131122}=\mp e^{2 X_3}(2b(X_{3})+b'(X_{3})),\nonumber\\
&(R\cdot P)_{131133}=-e^{4 X_3}(2b(X_{3})+b'(X_{3})),~(R\cdot P)_{131223}=-(R\cdot P)_{311223}=\mp 2,\nonumber\\
&(R\cdot P)_{132112}=-(R\cdot P)_{312112}=\pm \frac{1}{2}e^{2 X_3}(2b(X_{3})+b'(X_{3})),\nonumber\\
&(R\cdot P)_{132123}=-(R\cdot P)_{312123}=\pm 1,~(R\cdot P)_{133311}=e^{4 X_3}(2b(X_{3})+b'(X_{3})),\nonumber\\
&(R\cdot P)_{132323}=-(R\cdot P)_{312323}=\frac{1}{2}(2b(X_{3})+b'(X_{3})),~(R\cdot P)_{133322}=2b(X_{3})+b'(X_{3}),\nonumber\\
&(R\cdot P)_{133313}=-(R\cdot P)_{313313}=e^{2 X_3},~(R\cdot P)_{133333}=\pm e^{2 X_3}(2b(X_{3})+b'(X_{3})),\nonumber\\
&(R\cdot P)_{221223}=-(R\cdot P)_{222123}=\mp e^{-2 X_3}(2b(X_{3})+b'(X_{3})),\nonumber\\
&(R\cdot P)_{222312}=-(R\cdot P)_{223212}=\mp e^{-2 X_3}(2b(X_{3})+b'(X_{3})),\nonumber\\
&(R\cdot P)_{231211}=\pm e^{2 X_3}(2b(X_{3})+b'(X_{3})),~(R\cdot P)_{231213}=-(R\cdot P)_{321213}=\pm 1,\nonumber\\
&(R\cdot P)_{231222}=\pm e^{-2 X_3}(2b(X_{3})+b'(X_{3})),~(R\cdot P)_{231233}=2b(X_{3})+b'(X_{3}),\nonumber\\
&(R\cdot P)_{231223}=-(R\cdot P)_{321223}=-\frac{1}{2}(2b(X_{3})+b'(X_{3})),\nonumber\\
&(R\cdot P)_{232111}=\mp 2e^{2 X_3}(2b(X_{3})+b'(X_{3})),~(R\cdot P)_{232113}=-(R\cdot P)_{322113}=\mp 2,\nonumber\\
&(R\cdot P)_{232122}=\mp 2e^{-2 X_3}(2b(X_{3})+b'(X_{3})),~(R\cdot P)_{232133}=-2(2b(X_{3})+b'(X_{3})),\nonumber\\
&(R\cdot P)_{232212}=-(R\cdot P)_{322212}=\pm \frac{1}{2}e^{-2 X_3}(2b(X_{3})+b'(X_{3})),\nonumber\\
&(R\cdot P)_{232223}=-(R\cdot P)_{322223}=\mp e^{-4 X_3},~(R\cdot P)_{233323}=-(R\cdot P)_{323323}=e^{-2 X_3},\nonumber\\
&(R\cdot P)_{233312}=-(R\cdot P)_{323312}=\frac{1}{2}(2b(X_{3})+b'(X_{3})),\nonumber\\
&(R\cdot P)_{331223}=-(R\cdot P)_{332123}=(R\cdot P)_{332312}=-(R\cdot P)_{333212}=-(2b(X_{3})+b'(X_{3})).\nonumber
\end{align}
We notice that
\begin{align}
&(C\cdot R)_{111212}=\mp e^{4 X_3}(2b(X_{3})+b'(X_{3}))^2,~(C\cdot R)_{112323}=-e^{2 X_3}(2b(X_{3})+b'(X_{3}))^2,\\
&(C\cdot R)_{121312}=\mp 2e^{2 X_3}(2b(X_{3})+b'(X_{3})),~(C\cdot R)_{122311}=\pm 2e^{2 X_3}(2b(X_{3})+b'(X_{3})),\nonumber\\
&(C\cdot R)_{122322}=\pm 2e^{-2 X_3}(2b(X_{3})+b'(X_{3})),~(C\cdot R)_{122333}=2(2b(X_{3})+b'(X_{3})),\nonumber\\
&(C\cdot R)_{131212}=\mp 2e^{2 X_3}(2b(X_{3})+b'(X_{3})),~(C\cdot R)_{221212}=\mp (2b(X_{3})+b'(X_{3}))^2,\nonumber\\
&(C\cdot R)_{222323}=-e^{-2 X_3}(2b(X_{3})+b'(X_{3}))^2,~(C\cdot R)_{231211}=\pm 2e^{2 X_3}(2b(X_{3})+b'(X_{3})),\nonumber\\
&(C\cdot R)_{231222}=\pm 2e^{-2 X_3}(2b(X_{3})+b'(X_{3})),~(C\cdot R)_{231233}=2(2b(X_{3})+b'(X_{3})),\nonumber\\
&(C\cdot R)_{331212}=-e^{2 X_3}(2b(X_{3})+b'(X_{3}))^2,~(C\cdot R)_{332323}=\mp (2b(X_{3})+b'(X_{3}))^2;\nonumber
\end{align}
\begin{align}
&(C\cdot C)_{111212}=\mp e^{4 X_3}(2b(X_{3})+b'(X_{3}))^2,~(C\cdot C)_{112323}=-e^{2 X_3}(2b(X_{3})+b'(X_{3}))^2,\\
&(C\cdot C)_{121323}=e^{2 X_3}(2b(X_{3})+b'(X_{3}))^2,~(C\cdot C)_{131223}=-e^{2 X_3}(2b(X_{3})+b'(X_{3}))^2,\nonumber\\
&(C\cdot C)_{132312}=e^{2 X_3}(2b(X_{3})+b'(X_{3}))^2,~(C\cdot C)_{221212}=\mp (2b(X_{3})+b'(X_{3}))^2,\nonumber\\
&(C\cdot C)_{222323}=-e^{-2 X_3}(2b(X_{3})+b'(X_{3}))^2,~(C\cdot C)_{231312}=-e^{2 X_3}(2b(X_{3})+b'(X_{3}))^2;\nonumber\\
&(C\cdot K)_{111212}=\mp e^{4 X_3}(2b(X_{3})+b'(X_{3}))^2,~(C\cdot K)_{112323}=-e^{2 X_3}(2b(X_{3})+b'(X_{3}))^2,\\
&(C\cdot K)_{121323}=e^{2 X_3}(2b(X_{3})+b'(X_{3}))^2,~(C\cdot K)_{131223}=-e^{2 X_3}(2b(X_{3})+b'(X_{3}))^2,\nonumber\\
&(C\cdot K)_{132312}=e^{2 X_3}(2b(X_{3})+b'(X_{3}))^2,~(C\cdot K)_{221212}=\mp (2b(X_{3})+b'(X_{3}))^2,\nonumber\\
&(C\cdot K)_{222323}=-e^{-2 X_3}(2b(X_{3})+b'(X_{3}))^2,~(C\cdot K)_{231312}=-e^{2 X_3}(2b(X_{3})+b'(X_{3}))^2,\nonumber\\
&(C\cdot K)_{331212}=-e^{2 X_3}(2b(X_{3})+b'(X_{3}))^2,~(C\cdot K)_{332323}=\mp (2b(X_{3})+b'(X_{3}))^2;\nonumber\\
&(C\cdot W)_{111212}=\mp e^{4 X_3}(2b(X_{3})+b'(X_{3}))^2,~(C\cdot W)_{112323}=-e^{2 X_3}(2b(X_{3})+b'(X_{3}))^2,\\
&(C\cdot W)_{121312}=\mp 2e^{2 X_3}(2b(X_{3})+b'(X_{3})),~(C\cdot W)_{122311}=\pm 2e^{2 X_3}(2b(X_{3})+b'(X_{3})),\nonumber\\
&(C\cdot W)_{122322}=\pm 2e^{-2 X_3}(2b(X_{3})+b'(X_{3})),~(C\cdot W)_{122333}=2(2b(X_{3})+b'(X_{3})),\nonumber\\
&(C\cdot W)_{131212}=\mp 2e^{2 X_3}(2b(X_{3})+b'(X_{3})),~(C\cdot W)_{221212}=\mp (2b(X_{3})+b'(X_{3}))^2,\nonumber\\
&(C\cdot W)_{222323}=-e^{-2 X_3}(2b(X_{3})+b'(X_{3}))^2,~(C\cdot W)_{231211}=\pm 2e^{2 X_3}(2b(X_{3})+b'(X_{3})),\nonumber\\
&(C\cdot W)_{231222}=\pm 2e^{-2 X_3}(2b(X_{3})+b'(X_{3})),~(C\cdot W)_{231233}=2(2b(X_{3})+b'(X_{3})),\nonumber\\
&(C\cdot W)_{331212}=-e^{2 X_3}(2b(X_{3})+b'(X_{3}))^2,~(C\cdot W)_{332323}=\mp (2b(X_{3})+b'(X_{3}))^2;\nonumber\\
&(C\cdot P)_{111212}=-(C\cdot P)_{112112}=\mp e^{4 X_3}(2b(X_{3})+b'(X_{3}))^2,\\
&(C\cdot P)_{112323}=-(C\cdot P)_{113223}=-e^{2 X_3}(2b(X_{3})+b'(X_{3}))^2,\nonumber\\
&(C\cdot P)_{121112}=-(C\cdot P)_{211112}=\pm \frac{1}{2}e^{4 X_3}(2b(X_{3})+b'(X_{3}))^2,\nonumber\\
&(C\cdot P)_{121312}=-(C\cdot P)_{211312}=\mp 2e^{2 X_3}(2b(X_{3})+b'(X_{3})),\nonumber\\
&(C\cdot P)_{122212}=-(C\cdot P)_{212212}=\pm \frac{1}{2}(2b(X_{3})+b'(X_{3}))^2,\nonumber\\
&(C\cdot P)_{122311}=\pm 2e^{2 X_3}(2b(X_{3})+b'(X_{3})),~(C\cdot P)_{122322}=\pm 2e^{-2 X_3}(2b(X_{3})+b'(X_{3})),\nonumber\\
&(C\cdot P)_{122333}=2(2b(X_{3})+b'(X_{3})),~(C\cdot P)_{123233}=-(2b(X_{3})+b'(X_{3})),\nonumber\\
&(C\cdot P)_{123112}=-(C\cdot P)_{213112}=\pm e^{2 X_3}(2b(X_{3})+b'(X_{3})),\nonumber\\
&(C\cdot P)_{123123}=-(C\cdot P)_{213123}=-\frac{1}{2}e^{2 X_3}(2b(X_{3})+b'(X_{3}))^2,\nonumber
\end{align}
\begin{align}
&(C\cdot P)_{123211}=\mp e^{2 X_3}(2b(X_{3})+b'(X_{3})),~(C\cdot P)_{123222}=\mp e^{-2 X_3}(2b(X_{3})+b'(X_{3})),\nonumber\\
&(C\cdot P)_{131111}=\mp e^{6 X_3}(2b(X_{3})+b'(X_{3})),~(C\cdot P)_{131122}=\mp e^{2 X_3}(2b(X_{3})+b'(X_{3})),\nonumber\\
&(C\cdot P)_{131133}=-(C\cdot P)_{133311}=-e^{4 X_3}(2b(X_{3})+b'(X_{3})),\nonumber\\
&(C\cdot P)_{131212}=-(C\cdot P)_{311212}=\mp 2e^{2 X_3}(2b(X_{3})+b'(X_{3})),\nonumber\\
&(C\cdot P)_{132112}=-(C\cdot P)_{312112}=\pm e^{2 X_3}(2b(X_{3})+b'(X_{3})),\nonumber\\
&(C\cdot P)_{132123}=-(C\cdot P)_{312123}=\frac{1}{2}e^{2 X_3}(2b(X_{3})+b'(X_{3}))^2,\nonumber\\
&(C\cdot P)_{132312}=-(C\cdot P)_{312312}=\frac{1}{2}e^{2 X_3}(2b(X_{3})+b'(X_{3}))^2,\nonumber\\
&(C\cdot P)_{133322}=2b(X_{3})+b'(X_{3}),~(C\cdot P)_{133333}=\pm e^{2 X_3}(2b(X_{3})+b'(X_{3})),\nonumber\\
&(C\cdot P)_{221212}=-(C\cdot P)_{222112}=\mp (2b(X_{3})+b'(X_{3}))^2,\nonumber\\
&(C\cdot P)_{222323}=-(C\cdot P)_{223223}=-e^{-2 X_3}(2b(X_{3})+b'(X_{3}))^2,\nonumber\\
&(C\cdot P)_{231211}=\pm e^{2 X_3}(2b(X_{3})+b'(X_{3})),~(C\cdot P)_{231222}=\pm e^{-2 X_3}(2b(X_{3})+b'(X_{3})),\nonumber\\
&(C\cdot P)_{231233}=2b(X_{3})+b'(X_{3}),~(C\cdot P)_{232133}=-2(2b(X_{3})+b'(X_{3})),\nonumber\\
&(C\cdot P)_{231312}=-(C\cdot P)_{321312}=-\frac{1}{2}e^{2 X_3}(2b(X_{3})+b'(X_{3}))^2,\nonumber\\
&(C\cdot P)_{232111}=\mp 2e^{2 X_3}(2b(X_{3})+b'(X_{3})),~(C\cdot P)_{232122}=\mp 2e^{-2 X_3}(2b(X_{3})+b'(X_{3})),\nonumber\\
&(C\cdot P)_{232212}=-(C\cdot P)_{322212}=\mp e^{-2 X_3}(2b(X_{3})+b'(X_{3})),\nonumber\\
&(C\cdot P)_{232223}=-(C\cdot P)_{322223}=\frac{1}{2}e^{-2 X_3}(2b(X_{3})+b'(X_{3}))^2,\nonumber\\
&(C\cdot P)_{233312}=-(C\cdot P)_{323312}=2b(X_{3})+b'(X_{3}),\nonumber\\
&(C\cdot P)_{233323}=-(C\cdot P)_{323323}=\pm \frac{1}{2}(2b(X_{3})+b'(X_{3}))^2,\nonumber\\
&(C\cdot P)_{331212}=-(C\cdot P)_{332112}=-e^{2 X_3}(2b(X_{3})+b'(X_{3}))^2,\nonumber\\
&(C\cdot P)_{332323}=-(C\cdot P)_{333223}=\mp (2b(X_{3})+b'(X_{3}))^2.\nonumber
\end{align}
It is obvious that
\begin{align}
&(K\cdot R)_{111212}=\mp e^{4 X_3}(2b(X_{3})+b'(X_{3}))^2,~(K\cdot R)_{112323}=-e^{2 X_3}(2b(X_{3})+b'(X_{3}))^2,\\
&(K\cdot R)_{111223}=-(K\cdot R)_{112312}=\pm e^{2 X_3}(2b(X_{3})+b'(X_{3})),\nonumber\\
&(K\cdot R)_{121312}=(K\cdot R)_{131212}=\mp 2e^{2 X_3}(2b(X_{3})+b'(X_{3})),\nonumber\\
&(K\cdot R)_{121323}=(K\cdot R)_{131223}=\pm 2,\nonumber\\
&(K\cdot R)_{122311}=(K\cdot R)_{231211}=\pm 2e^{2 X_3}(2b(X_{3})+b'(X_{3})),\nonumber
\end{align}
\begin{align}
&(K\cdot R)_{122313}=(K\cdot R)_{231213}=\mp 2,\nonumber\\
&(K\cdot R)_{122322}=(K\cdot R)_{231222}=\pm 2e^{-2 X_3}(2b(X_{3})+b'(X_{3})),\nonumber\\
&(K\cdot R)_{122333}=(K\cdot R)_{231233}=2(2b(X_{3})+b'(X_{3})),\nonumber\\
&(K\cdot R)_{221212}=\mp (2b(X_{3})+b'(X_{3}))^2,~(K\cdot R)_{331212}=-e^{2 X_3}(2b(X_{3})+b'(X_{3}))^2,\nonumber\\
&(K\cdot R)_{221223}=-(K\cdot R)_{222312}=\pm e^{-2 X_3}(2b(X_{3})+b'(X_{3})),\nonumber\\
&(K\cdot R)_{222323}=-e^{-2 X_3}(2b(X_{3})+b'(X_{3}))^2,~(K\cdot R)_{332323}=\mp (2b(X_{3})+b'(X_{3}))^2,\nonumber\\
&(K\cdot R)_{331223}=-(K\cdot R)_{332312}=2b(X_{3})+b'(X_{3});\nonumber\\
&(K\cdot C)_{111212}=\mp e^{4 X_3}(2b(X_{3})+b'(X_{3}))^2,~(K\cdot C)_{112323}=-e^{2 X_3}(2b(X_{3})+b'(X_{3}))^2,\\
&(K\cdot C)_{111223}=-(K\cdot C)_{112312}=\pm e^{2 X_3}(2b(X_{3})+b'(X_{3})),\nonumber\\
&(K\cdot C)_{121312}=-(K\cdot C)_{131212}=\pm e^{2 X_3}(2b(X_{3})+b'(X_{3})),\nonumber\\
&(K\cdot C)_{121323}=-(K\cdot C)_{131223}=e^{2 X_3}(2b(X_{3})+b'(X_{3}))^2,\nonumber\\
&(K\cdot C)_{132312}=-(K\cdot C)_{231312}=e^{2 X_3}(2b(X_{3})+b'(X_{3}))^2,\nonumber\\
&(K\cdot C)_{132323}=-(K\cdot C)_{231323}=-(2b(X_{3})+b'(X_{3})),\nonumber\\
&(K\cdot C)_{221212}=\mp (2b(X_{3})+b'(X_{3}))^2,~(K\cdot C)_{331212}=-e^{2 X_3}(2b(X_{3})+b'(X_{3}))^2,\nonumber\\
&(K\cdot C)_{221223}=-(K\cdot C)_{222312}=\pm e^{-2 X_3}(2b(X_{3})+b'(X_{3})),\nonumber\\
&(K\cdot C)_{222323}=-e^{-2 X_3}(2b(X_{3})+b'(X_{3}))^2,~(K\cdot C)_{332323}=\mp (2b(X_{3})+b'(X_{3}))^2,\nonumber\\
&(K\cdot C)_{331223}=-(K\cdot C)_{332312}=2b(X_{3})+b'(X_{3});\nonumber\\
&(K\cdot K)_{111212}=\mp e^{4 X_3}(2b(X_{3})+b'(X_{3}))^2,~(K\cdot K)_{112323}=-e^{2 X_3}(2b(X_{3})+b'(X_{3}))^2,\\
&(K\cdot K)_{111223}=-(K\cdot K)_{112312}=\pm e^{2 X_3}(2b(X_{3})+b'(X_{3})),\nonumber\\
&(K\cdot K)_{121312}=-(K\cdot K)_{131212}=\pm e^{2 X_3}(2b(X_{3})+b'(X_{3})),\nonumber\\
&(K\cdot K)_{121323}=-(K\cdot K)_{131223}=e^{2 X_3}(2b(X_{3})+b'(X_{3}))^2,\nonumber\\
&(K\cdot K)_{132312}=-(K\cdot K)_{231312}=e^{2 X_3}(2b(X_{3})+b'(X_{3}))^2,\nonumber\\
&(K\cdot K)_{132323}=-(K\cdot K)_{231323}=-(2b(X_{3})+b'(X_{3})),\nonumber\\
&(K\cdot K)_{221212}=\mp (2b(X_{3})+b'(X_{3}))^2,~(K\cdot K)_{331212}=-e^{2 X_3}(2b(X_{3})+b'(X_{3}))^2,\nonumber\\
&(K\cdot K)_{221223}=-(K\cdot K)_{222312}=\pm e^{-2 X_3}(2b(X_{3})+b'(X_{3})),\nonumber\\
&(K\cdot K)_{222323}=-e^{-2 X_3}(2b(X_{3})+b'(X_{3}))^2,~(K\cdot K)_{332323}=\mp (2b(X_{3})+b'(X_{3}))^2,\nonumber\\
&(K\cdot K)_{331223}=-(K\cdot K)_{332312}=2b(X_{3})+b'(X_{3});\nonumber
\end{align}
\begin{align}
&(K\cdot W)_{111212}=\mp e^{4 X_3}(2b(X_{3})+b'(X_{3}))^2,~(K\cdot W)_{112323}=-e^{2 X_3}(2b(X_{3})+b'(X_{3}))^2,\\
&(K\cdot W)_{111223}=-(K\cdot W)_{112312}=\pm e^{2 X_3}(2b(X_{3})+b'(X_{3})),\nonumber\\
&(K\cdot W)_{121312}=(K\cdot W)_{131212}=\mp 2e^{2 X_3}(2b(X_{3})+b'(X_{3})),\nonumber\\
&(K\cdot W)_{121323}=(K\cdot W)_{131223}=\pm 2,\nonumber\\
&(K\cdot W)_{122311}=(K\cdot W)_{231211}=\pm 2e^{2 X_3}(2b(X_{3})+b'(X_{3})),\nonumber\\
&(K\cdot W)_{122313}=(K\cdot W)_{231213}=\mp 2,\nonumber\\
&(K\cdot W)_{122322}=(K\cdot W)_{231222}=\pm 2e^{-2 X_3}(2b(X_{3})+b'(X_{3})),\nonumber\\
&(K\cdot W)_{122333}=(K\cdot W)_{231233}=2(2b(X_{3})+b'(X_{3})),\nonumber\\
&(K\cdot W)_{221212}=\mp (2b(X_{3})+b'(X_{3}))^2,~(K\cdot W)_{331212}=-e^{2 X_3}(2b(X_{3})+b'(X_{3}))^2,\nonumber\\
&(K\cdot W)_{221223}=-(K\cdot W)_{222312}=\pm e^{-2 X_3}(2b(X_{3})+b'(X_{3})),\nonumber\\
&(K\cdot W)_{222323}=-e^{-2 X_3}(2b(X_{3})+b'(X_{3}))^2,~(K\cdot W)_{332323}=\mp (2b(X_{3})+b'(X_{3}))^2,\nonumber\\
&(K\cdot W)_{331223}=-(K\cdot W)_{332312}=2b(X_{3})+b'(X_{3});\nonumber\\
&(K\cdot P)_{111212}=-(K\cdot P)_{112112}=\mp e^{4 X_3}(2b(X_{3})+b'(X_{3}))^2,\\
&(K\cdot P)_{111223}=-(K\cdot P)_{112123}=\pm e^{2 X_3}(2b(X_{3})+b'(X_{3})),\nonumber\\
&(K\cdot P)_{112312}=-(K\cdot P)_{113212}=\mp e^{2 X_3}(2b(X_{3})+b'(X_{3})),\nonumber\\
&(K\cdot P)_{112323}=-(K\cdot P)_{113223}=-e^{2 X_3}(2b(X_{3})+b'(X_{3}))^2,\nonumber\\
&(K\cdot P)_{121112}=-(K\cdot P)_{211112}=\pm \frac{1}{2}e^{4 X_3}(2b(X_{3})+b'(X_{3}))^2,\nonumber\\
&(K\cdot P)_{121123}=-(K\cdot P)_{211123}=\mp \frac{1}{2}e^{2 X_3}(2b(X_{3})+b'(X_{3})),\nonumber\\
&(K\cdot P)_{121312}=-(K\cdot P)_{211312}=\mp 2e^{2 X_3}(2b(X_{3})+b'(X_{3})),\nonumber\\
&(K\cdot P)_{121323}=-(K\cdot P)_{211323}=\pm 2,\nonumber\\
&(K\cdot P)_{122212}=-(K\cdot P)_{212212}=\pm \frac{1}{2}(2b(X_{3})+b'(X_{3}))^2,\nonumber\\
&(K\cdot P)_{122223}=-(K\cdot P)_{212223}=\mp \frac{1}{2}e^{-2 X_3}(2b(X_{3})+b'(X_{3})),\nonumber\\
&(K\cdot P)_{122311}=\pm 2e^{2 X_3}(2b(X_{3})+b'(X_{3})),~(K\cdot P)_{123211}=\mp e^{2 X_3}(2b(X_{3})+b'(X_{3})),\nonumber\\
&(K\cdot P)_{122313}=-(K\cdot P)_{212313}=\mp 2,~(K\cdot P)_{123213}=-(K\cdot P)_{213213}=\pm 1,\nonumber\\
&(K\cdot P)_{122322}=\pm 2e^{-2 X_3}(2b(X_{3})+b'(X_{3})),~(K\cdot P)_{123222}=\mp e^{-2 X_3}(2b(X_{3})+b'(X_{3})),\nonumber\\
&(K\cdot P)_{122333}=2(2b(X_{3})+b'(X_{3})),~(K\cdot P)_{123233}=-(2b(X_{3})+b'(X_{3})),\nonumber
\end{align}
\begin{align}
&(K\cdot P)_{123112}=-(K\cdot P)_{213112}=\pm \frac{1}{2}e^{2 X_3}(2b(X_{3})+b'(X_{3})),\nonumber\\
&(K\cdot P)_{123123}=-(K\cdot P)_{213123}=\mp 1-\frac{1}{2}e^{2 X_3}(2b(X_{3})+b'(X_{3}))^2,\nonumber\\
&(K\cdot P)_{131111}=\mp e^{6 X_3}(2b(X_{3})+b'(X_{3})),~(K\cdot P)_{133311}=e^{4 X_3}(2b(X_{3})+b'(X_{3})),\nonumber\\
&(K\cdot P)_{131113}=-(K\cdot P)_{311113}=\pm e^{4 X_3},~(K\cdot P)_{133313}=-(K\cdot P)_{313313}=-e^{2 X_3},\nonumber\\
&(K\cdot P)_{131122}=\mp e^{2 X_3}(2b(X_{3})+b'(X_{3})),~(K\cdot P)_{133322}=2b(X_{3})+b'(X_{3}),\nonumber\\
&(K\cdot P)_{131133}=-e^{4 X_3}(2b(X_{3})+b'(X_{3})),~(K\cdot P)_{133333}=\pm e^{2 X_3}(2b(X_{3})+b'(X_{3})),\nonumber\\
&(K\cdot P)_{131212}=-(K\cdot P)_{311212}=\mp 2e^{2 X_3}(2b(X_{3})+b'(X_{3})),\nonumber\\
&(K\cdot P)_{131223}=-(K\cdot P)_{311223}=\pm 2,\nonumber\\
&(K\cdot P)_{132112}=-(K\cdot P)_{312112}=\pm \frac{3}{2}e^{2 X_3}(2b(X_{3})+b'(X_{3})),\nonumber\\
&(K\cdot P)_{132123}=-(K\cdot P)_{312123}=\mp 1+\frac{1}{2}e^{2 X_3}(2b(X_{3})+b'(X_{3}))^2,\nonumber\\
&(K\cdot P)_{132312}=-(K\cdot P)_{312312}=\frac{1}{2}e^{2 X_3}(2b(X_{3})+b'(X_{3}))^2,\nonumber\\
&(K\cdot P)_{132323}=-(K\cdot P)_{312323}=-\frac{1}{2}(2b(X_{3})+b'(X_{3})),\nonumber\\
&(K\cdot P)_{221212}=-(K\cdot P)_{222112}=\mp (2b(X_{3})+b'(X_{3}))^2,\nonumber\\
&(K\cdot P)_{221223}=-(K\cdot P)_{222123}=\pm e^{-2 X_3}(2b(X_{3})+b'(X_{3})),\nonumber\\
&(K\cdot P)_{222312}=-(K\cdot P)_{223212}=\mp e^{-2 X_3}(2b(X_{3})+b'(X_{3})),\nonumber\\
&(K\cdot P)_{222323}=-(K\cdot P)_{223223}=-e^{-2 X_3}(2b(X_{3})+b'(X_{3}))^2,\nonumber\\
&(K\cdot P)_{231211}=\pm e^{2 X_3}(2b(X_{3})+b'(X_{3})),~(K\cdot P)_{232111}=\mp 2e^{2 X_3}(2b(X_{3})+b'(X_{3})),\nonumber\\
&(K\cdot P)_{231213}=-(K\cdot P)_{321213}=\mp 1,~(K\cdot P)_{232113}=-(K\cdot P)_{322113}=\pm 2,\nonumber\\
&(K\cdot P)_{231222}=\pm e^{-2 X_3}(2b(X_{3})+b'(X_{3})),~(K\cdot P)_{232122}=\mp 2e^{-2 X_3}(2b(X_{3})+b'(X_{3})),\nonumber\\
&(K\cdot P)_{231233}=2b(X_{3})+b'(X_{3}),~(K\cdot P)_{232133}=-2(2b(X_{3})+b'(X_{3})),\nonumber\\
&(K\cdot P)_{231312}=-(K\cdot P)_{321312}=-\frac{1}{2}e^{2 X_3}(2b(X_{3})+b'(X_{3}))^2,\nonumber\\
&(K\cdot P)_{231323}=-(K\cdot P)_{321323}=\frac{1}{2}(2b(X_{3})+b'(X_{3})),\nonumber\\
&(K\cdot P)_{232212}=-(K\cdot P)_{322212}=\mp \frac{1}{2}e^{-2 X_3}(2b(X_{3})+b'(X_{3})),\nonumber\\
&(K\cdot P)_{232223}=-(K\cdot P)_{322223}=\pm e^{-4 X_3}+\frac{1}{2}e^{-2 X_3}(2b(X_{3})+b'(X_{3}))^2,\nonumber\\
&(K\cdot P)_{233312}=-(K\cdot P)_{323312}=\frac{3}{2}(2b(X_{3})+b'(X_{3})),\nonumber
\end{align}
\begin{align}
&(K\cdot P)_{233323}=-(K\cdot P)_{323323}=-e^{-2 X_3}\pm \frac{1}{2}(2b(X_{3})+b'(X_{3}))^2,\nonumber\\
&(K\cdot P)_{331212}=-(K\cdot P)_{332112}=-e^{2 X_3}(2b(X_{3})+b'(X_{3}))^2,\nonumber\\
&(K\cdot P)_{331223}=-(K\cdot P)_{332123}=2b(X_{3})+b'(X_{3}),\nonumber\\
&(K\cdot P)_{332312}=-(K\cdot P)_{333212}=-(2b(X_{3})+b'(X_{3})),\nonumber\\
&(K\cdot P)_{332323}=-(K\cdot P)_{333223}=\mp (2b(X_{3})+b'(X_{3}))^2.\nonumber
\end{align}
Computations show that
\begin{align}
&(W\cdot R)_{111223}=\mp \frac{2}{3}e^{2 X_3}(2b(X_{3})+b'(X_{3})),~(W\cdot R)_{112312}=\mp \frac{4}{3}e^{2 X_3}(2b(X_{3})+b'(X_{3})),\\
&(W\cdot R)_{121323}=(W\cdot R)_{131223}=\mp \frac{4}{3},~(W\cdot R)_{122313}=(W\cdot R)_{231213}=\pm \frac{4}{3},\nonumber\\
&(W\cdot R)_{122311}=(W\cdot R)_{231211}=\pm 2e^{2 X_3}(2b(X_{3})+b'(X_{3})),\nonumber\\
&(W\cdot R)_{122322}=(W\cdot R)_{231222}=\pm 2e^{-2 X_3}(2b(X_{3})+b'(X_{3})),\nonumber\\
&(W\cdot R)_{122333}=(W\cdot R)_{231233}=2(2b(X_{3})+b'(X_{3})),\nonumber\\
&(W\cdot R)_{221223}=\mp \frac{2}{3}e^{-2 X_3}(2b(X_{3})+b'(X_{3})),~(W\cdot R)_{222312}=\mp \frac{4}{3}e^{-2 X_3}(2b(X_{3})+b'(X_{3})),\nonumber\\
&(W\cdot R)_{331223}=-\frac{2}{3}(2b(X_{3})+b'(X_{3})),~(W\cdot R)_{332312}=-\frac{4}{3}(2b(X_{3})+b'(X_{3}));\nonumber\\
&(W\cdot C)_{111223}=\mp \frac{2}{3}e^{2 X_3}(2b(X_{3})+b'(X_{3})),~(W\cdot C)_{112312}=\mp \frac{4}{3}e^{2 X_3}(2b(X_{3})+b'(X_{3})),\\
&(W\cdot C)_{121312}=-(W\cdot C)_{131212}=\pm \frac{4}{3}e^{2 X_3}(2b(X_{3})+b'(X_{3})),\nonumber\\
&(W\cdot C)_{132323}=-(W\cdot C)_{231323}=\frac{2}{3}(2b(X_{3})+b'(X_{3})),\nonumber\\
&(W\cdot C)_{221223}=\mp \frac{2}{3}e^{-2 X_3}(2b(X_{3})+b'(X_{3})),~(W\cdot C)_{222312}=\mp \frac{4}{3}e^{-2 X_3}(2b(X_{3})+b'(X_{3})),\nonumber\\
&(W\cdot C)_{331223}=-\frac{2}{3}(2b(X_{3})+b'(X_{3})),~(W\cdot C)_{332312}=-\frac{4}{3}(2b(X_{3})+b'(X_{3}));\nonumber\\
&(W\cdot K)_{111223}=\mp \frac{2}{3}e^{2 X_3}(2b(X_{3})+b'(X_{3})),~(W\cdot K)_{112312}=\mp \frac{4}{3}e^{2 X_3}(2b(X_{3})+b'(X_{3})),\\
&(W\cdot K)_{121312}=-(W\cdot K)_{131212}=\pm \frac{4}{3}e^{2 X_3}(2b(X_{3})+b'(X_{3})),\nonumber\\
&(W\cdot K)_{132323}=-(W\cdot K)_{231323}=\frac{2}{3}(2b(X_{3})+b'(X_{3})),\nonumber
\end{align}
\begin{align}
&(W\cdot K)_{221223}=\mp \frac{2}{3}e^{-2 X_3}(2b(X_{3})+b'(X_{3})),~(W\cdot K)_{222312}=\mp \frac{4}{3}e^{-2 X_3}(2b(X_{3})+b'(X_{3})),\nonumber\\
&(W\cdot K)_{331223}=-\frac{2}{3}(2b(X_{3})+b'(X_{3})),~(W\cdot K)_{332312}=-\frac{4}{3}(2b(X_{3})+b'(X_{3}));\nonumber\\
&(W\cdot W)_{111223}=\mp \frac{2}{3}e^{2 X_3}(2b(X_{3})+b'(X_{3})),~(W\cdot W)_{112312}=\mp \frac{4}{3}e^{2 X_3}(2b(X_{3})+b'(X_{3})),\\
&(W\cdot W)_{121323}=(W\cdot W)_{131223}=\mp \frac{4}{3},~(W\cdot W)_{122313}=(W\cdot W)_{231213}=\pm \frac{4}{3},\nonumber\\
&(W\cdot W)_{122311}=(W\cdot W)_{231211}=\pm 2e^{2 X_3}(2b(X_{3})+b'(X_{3})),\nonumber\\
&(W\cdot W)_{122322}=(W\cdot W)_{231222}=\pm 2e^{-2 X_3}(2b(X_{3})+b'(X_{3})),\nonumber\\
&(W\cdot W)_{122333}=(W\cdot W)_{231233}=2(2b(X_{3})+b'(X_{3})),\nonumber\\
&(W\cdot W)_{221223}=\mp \frac{2}{3}e^{-2 X_3}(2b(X_{3})+b'(X_{3})),~(W\cdot W)_{222312}=\mp \frac{4}{3}e^{-2 X_3}(2b(X_{3})+b'(X_{3})),\nonumber\\
&(W\cdot W)_{331223}=-\frac{2}{3}(2b(X_{3})+b'(X_{3})),~(W\cdot W)_{332312}=-\frac{4}{3}(2b(X_{3})+b'(X_{3}));\nonumber\\
&(W\cdot P)_{111223}=-(W\cdot P)_{112123}=\mp \frac{2}{3}e^{2 X_3}(2b(X_{3})+b'(X_{3})),\\
&(W\cdot P)_{112312}=-(W\cdot P)_{113212}=\mp \frac{4}{3}e^{2 X_3}(2b(X_{3})+b'(X_{3})),\nonumber\\
&(W\cdot P)_{121123}=\pm \frac{1}{3}e^{2 X_3}(2b(X_{3})+b'(X_{3})),~(W\cdot P)_{122223}=\pm \frac{1}{3}e^{-2 X_3}(2b(X_{3})+b'(X_{3})),\nonumber\\
&(W\cdot P)_{121323}=-(W\cdot P)_{122313}=(W\cdot P)_{131223}=(W\cdot P)_{232113}=\mp \frac{4}{3},\nonumber\\
&(W\cdot P)_{123123}=-(W\cdot P)_{123213}=(W\cdot P)_{132123}=(W\cdot P)_{231213}=\pm \frac{2}{3},\nonumber\\
&(W\cdot P)_{122311}=\pm 2e^{2 X_3}(2b(X_{3})+b'(X_{3})),~(W\cdot P)_{123211}=\mp e^{2 X_3}(2b(X_{3})+b'(X_{3})),\nonumber\\
&(W\cdot P)_{122322}=\pm 2e^{-2 X_3}(2b(X_{3})+b'(X_{3})),~(W\cdot P)_{123222}=\mp e^{-2 X_3}(2b(X_{3})+b'(X_{3})),\nonumber\\
&(W\cdot P)_{122333}=2(2b(X_{3})+b'(X_{3})),~(W\cdot P)_{123233}=-(2b(X_{3})+b'(X_{3})),\nonumber\\
&(W\cdot P)_{123112}=\mp \frac{2}{3}e^{2 X_3}(2b(X_{3})+b'(X_{3})),~(W\cdot P)_{131111}=\mp e^{6 X_3}(2b(X_{3})+b'(X_{3})),\nonumber\\
&(W\cdot P)_{131113}=\mp \frac{2}{3}e^{4 X_3},~(W\cdot P)_{133313}=\frac{2}{3}e^{2 X_3},~(W\cdot P)_{131122}=\mp e^{2 X_3}(2b(X_{3})+b'(X_{3})),\nonumber\\
&(W\cdot P)_{131133}=-e^{4 X_3}(2b(X_{3})+b'(X_{3})),~(W\cdot P)_{132112}=\pm \frac{2}{3}e^{2 X_3}(2b(X_{3})+b'(X_{3})),\nonumber\\
&(W\cdot P)_{132323}=\frac{1}{3}(2b(X_{3})+b'(X_{3})),~(W\cdot P)_{133311}=e^{4 X_3}(2b(X_{3})+b'(X_{3})),\nonumber
\end{align}
\begin{align}
&(W\cdot P)_{133322}=2b(X_{3})+b'(X_{3}),~(W\cdot P)_{133333}=\pm e^{2 X_3}(2b(X_{3})+b'(X_{3})),\nonumber\\
&(W\cdot P)_{221223}=-(W\cdot P)_{222123}=\mp \frac{2}{3}e^{-2 X_3}(2b(X_{3})+b'(X_{3})),\nonumber\\
&(W\cdot P)_{222312}=-(W\cdot P)_{223212}=\mp \frac{4}{3}e^{-2 X_3}(2b(X_{3})+b'(X_{3})),\nonumber\\
&(W\cdot P)_{231211}=\pm e^{2 X_3}(2b(X_{3})+b'(X_{3})),~(W\cdot P)_{232111}=\mp 2e^{2 X_3}(2b(X_{3})+b'(X_{3})),\nonumber\\
&(W\cdot P)_{231222}=\pm e^{-2 X_3}(2b(X_{3})+b'(X_{3})),~(W\cdot P)_{232122}=\mp 2e^{-2 X_3}(2b(X_{3})+b'(X_{3})),\nonumber\\
&(W\cdot P)_{231233}=2b(X_{3})+b'(X_{3}),~(W\cdot P)_{232133}=-2(2b(X_{3})+b'(X_{3})),\nonumber\\
&(W\cdot P)_{231323}=-\frac{1}{3}(2b(X_{3})+b'(X_{3})),~(W\cdot P)_{232212}=\pm \frac{2}{3}e^{-2 X_3}(2b(X_{3})+b'(X_{3})),\nonumber\\
&(W\cdot P)_{232223}=\mp \frac{2}{3}e^{-4 X_3},~(W\cdot P)_{233323}=\frac{2}{3}e^{-2 X_3},~(W\cdot P)_{233312}=\frac{2}{3}(2b(X_{3})+b'(X_{3})),\nonumber\\
&(W\cdot P)_{331223}=-(W\cdot P)_{332123}=-\frac{2}{3}(2b(X_{3})+b'(X_{3})),\nonumber\\
&(W\cdot P)_{332312}=-(W\cdot P)_{333212}=-\frac{4}{3}(2b(X_{3})+b'(X_{3})).\nonumber
\end{align}
We find that
\begin{align}
&(P\cdot R)_{111223}=(P\cdot R)_{112312}=\mp e^{2 X_3}(2b(X_{3})+b'(X_{3})),\\
&(P\cdot R)_{111313}=2e^{6 X_3}(2b(X_{3})+b'(X_{3}))^2,\nonumber\\
&(P\cdot R)_{121213}=-(P\cdot R)_{211213}=\mp e^{2 X_3}(2b(X_{3})+b'(X_{3})),\nonumber\\
&(P\cdot R)_{121312}=-(P\cdot R)_{211312}=\mp \frac{1}{2}e^{2 X_3}(2b(X_{3})+b'(X_{3})),\nonumber\\
&(P\cdot R)_{121323}=\mp 1+\frac{1}{2}e^{2 X_3}(2b(X_{3})+b'(X_{3}))^2,\nonumber\\
&(P\cdot R)_{122311}=-(P\cdot R)_{212311}=\pm 2e^{2 X_3}(2b(X_{3})+b'(X_{3})),\nonumber\\
&(P\cdot R)_{122313}=-(P\cdot R)_{212313}=\pm 1,\nonumber\\
&(P\cdot R)_{122322}=-(P\cdot R)_{212322}=\pm 2e^{-2 X_3}(2b(X_{3})+b'(X_{3})),\nonumber\\
&(P\cdot R)_{122333}=-(P\cdot R)_{212333}=2(2b(X_{3})+b'(X_{3})),\nonumber\\
&(P\cdot R)_{131212}=-(P\cdot R)_{311212}=\mp \frac{1}{2}e^{2 X_3}(2b(X_{3})+b'(X_{3})),\nonumber\\
&(P\cdot R)_{131223}=-(P\cdot R)_{311223}=\mp 1,\nonumber\\
&(P\cdot R)_{131313}=e^{4 X_3}(2b(X_{3})+b'(X_{3})),\nonumber\\
&(P\cdot R)_{132323}=-(P\cdot R)_{312323}=\frac{1}{2}(2b(X_{3})+b'(X_{3})),\nonumber\\
&(P\cdot R)_{211323}=\pm 1+\frac{1}{2}e^{2 X_3}(2b(X_{3})+b'(X_{3}))^2,\nonumber
\end{align}
\begin{align}
&(P\cdot R)_{221223}=(P\cdot R)_{222312}=\mp e^{-2 X_3}(2b(X_{3})+b'(X_{3})),\nonumber\\
&(P\cdot R)_{221313}=e^{2 X_3}(2b(X_{3})+b'(X_{3}))^2,\nonumber\\
&(P\cdot R)_{231211}=-(P\cdot R)_{321211}=\pm 2e^{2 X_3}(2b(X_{3})+b'(X_{3})),\nonumber\\
&(P\cdot R)_{231213}=-(P\cdot R)_{321213}=\pm 1,\nonumber\\
&(P\cdot R)_{231222}=-(P\cdot R)_{321222}=\pm 2e^{-2 X_3}(2b(X_{3})+b'(X_{3})),\nonumber\\
&(P\cdot R)_{231233}=-(P\cdot R)_{321233}=2(2b(X_{3})+b'(X_{3})),\nonumber\\
&(P\cdot R)_{231312}=(P\cdot R)_{321312}=\frac{1}{2}e^{2 X_3}(2b(X_{3})+b'(X_{3}))^2,\nonumber\\
&(P\cdot R)_{231323}=-\frac{1}{2}(2b(X_{3})+b'(X_{3})),\nonumber\\
&(P\cdot R)_{232313}=-(P\cdot R)_{322313}=2b(X_{3})+b'(X_{3}),\nonumber\\
&(P\cdot R)_{311313}=-3e^{4 X_3}(2b(X_{3})+b'(X_{3})),\nonumber\\
&(P\cdot R)_{321323}=-\frac{3}{2}(2b(X_{3})+b'(X_{3})),\nonumber\\
&(P\cdot R)_{331223}=(P\cdot R)_{332312}=-(2b(X_{3})+b'(X_{3})),\nonumber\\
&(P\cdot R)_{331313}=\pm 2e^{4 X_3}(2b(X_{3})+b'(X_{3}))^2;\nonumber\\
&(P\cdot C)_{111223}=(P\cdot C)_{112312}=\mp e^{2 X_3}(2b(X_{3})+b'(X_{3})),\\
&(P\cdot C)_{111313}=2e^{6 X_3}(2b(X_{3})+b'(X_{3}))^2,\nonumber\\
&(P\cdot C)_{121312}=-(P\cdot C)_{211312}=\pm e^{2 X_3}(2b(X_{3})+b'(X_{3})),\nonumber\\
&(P\cdot C)_{121323}=e^{2 X_3}(2b(X_{3})+b'(X_{3}))^2,\nonumber\\
&(P\cdot C)_{122313}=-(P\cdot C)_{212313}=e^{2 X_3}(2b(X_{3})+b'(X_{3}))^2,\nonumber\\
&(P\cdot C)_{131212}=-(P\cdot C)_{311212}=\mp e^{2 X_3}(2b(X_{3})+b'(X_{3})),\nonumber\\
&(P\cdot C)_{131223}=-(P\cdot C)_{311223}=-\frac{1}{2}e^{2 X_3}(2b(X_{3})+b'(X_{3}))^2,\nonumber\\
&(P\cdot C)_{131313}=(P\cdot C)_{311313}=-e^{4 X_3}(2b(X_{3})+b'(X_{3})),\nonumber\\
&(P\cdot C)_{132312}=-(P\cdot C)_{312312}=\frac{1}{2}e^{2 X_3}(2b(X_{3})+b'(X_{3}))^2,\nonumber\\
&(P\cdot C)_{221223}=(P\cdot C)_{222312}=\mp e^{-2 X_3}(2b(X_{3})+b'(X_{3})),\nonumber\\
&(P\cdot C)_{221313}=e^{2 X_3}(2b(X_{3})+b'(X_{3}))^2,\nonumber\\
&(P\cdot C)_{231213}=-(P\cdot C)_{321213}=-e^{2 X_3}(2b(X_{3})+b'(X_{3}))^2,\nonumber\\
&(P\cdot C)_{231323}=(P\cdot C)_{321323}=-(2b(X_{3})+b'(X_{3})),\nonumber\\
&(P\cdot C)_{321312}=e^{2 X_3}(2b(X_{3})+b'(X_{3}))^2,\nonumber\\
&(P\cdot C)_{331223}=(P\cdot C)_{332312}=-(2b(X_{3})+b'(X_{3})),\nonumber\\
&(P\cdot C)_{331313}=\pm 2e^{4 X_3}(2b(X_{3})+b'(X_{3}))^2;\nonumber
\end{align}
\begin{align}
&(P\cdot K)_{111223}=(P\cdot K)_{112312}=\mp e^{2 X_3}(2b(X_{3})+b'(X_{3})),\\
&(P\cdot K)_{111313}=2e^{6 X_3}(2b(X_{3})+b'(X_{3}))^2,\nonumber\\
&(P\cdot K)_{121213}=-(P\cdot K)_{211213}=\mp e^{2 X_3}(2b(X_{3})+b'(X_{3})),\nonumber\\
&(P\cdot K)_{121312}=-(P\cdot K)_{211312}=\pm \frac{1}{2}e^{2 X_3}(2b(X_{3})+b'(X_{3})),\nonumber\\
&(P\cdot K)_{121323}=\pm 1+e^{2 X_3}(2b(X_{3})+b'(X_{3}))^2,~(P\cdot K)_{211323}=\mp 1,\nonumber\\
&(P\cdot K)_{122313}=-(P\cdot K)_{212313}=\mp 1+e^{2 X_3}(2b(X_{3})+b'(X_{3}))^2,\nonumber\\
&(P\cdot K)_{131212}=-(P\cdot K)_{311212}=\mp \frac{3}{2}e^{2 X_3}(2b(X_{3})+b'(X_{3})),\nonumber\\
&(P\cdot K)_{131223}=-(P\cdot K)_{311223}=\pm 1-\frac{1}{2}e^{2 X_3}(2b(X_{3})+b'(X_{3}))^2,\nonumber\\
&(P\cdot K)_{131313}=-3e^{4 X_3}(2b(X_{3})+b'(X_{3})),\nonumber\\
&(P\cdot K)_{132312}=-(P\cdot K)_{312312}=\frac{1}{2}e^{2 X_3}(2b(X_{3})+b'(X_{3}))^2,\nonumber\\
&(P\cdot K)_{132323}=-(P\cdot K)_{312323}=-\frac{1}{2}(2b(X_{3})+b'(X_{3})),\nonumber\\
&(P\cdot K)_{221223}=(P\cdot K)_{222312}=\mp e^{-2 X_3}(2b(X_{3})+b'(X_{3})),\nonumber\\
&(P\cdot K)_{221313}=e^{2 X_3}(2b(X_{3})+b'(X_{3}))^2,\nonumber\\
&(P\cdot K)_{231213}=-(P\cdot K)_{321213}=\mp 1-e^{2 X_3}(2b(X_{3})+b'(X_{3}))^2,\nonumber\\
&(P\cdot K)_{231323}=-\frac{3}{2}(2b(X_{3})+b'(X_{3})),\nonumber\\
&(P\cdot K)_{232313}=-(P\cdot K)_{322313}=-(2b(X_{3})+b'(X_{3})),\nonumber\\
&(P\cdot K)_{311313}=e^{4 X_3}(2b(X_{3})+b'(X_{3})),\nonumber\\
&(P\cdot K)_{321312}=e^{2 X_3}(2b(X_{3})+b'(X_{3}))^2,\nonumber\\
&(P\cdot K)_{321323}=-\frac{1}{2}(2b(X_{3})+b'(X_{3})),\nonumber\\
&(P\cdot K)_{331223}=(P\cdot K)_{332312}=-(2b(X_{3})+b'(X_{3})),\nonumber\\
&(P\cdot K)_{331313}=\pm 2e^{4 X_3}(2b(X_{3})+b'(X_{3}))^2;\nonumber\\
&(P\cdot W)_{111223}=(P\cdot W)_{112312}=\mp e^{2 X_3}(2b(X_{3})+b'(X_{3})),\\
&(P\cdot W)_{111313}=2e^{6 X_3}(2b(X_{3})+b'(X_{3}))^2,\nonumber\\
&(P\cdot W)_{121213}=-(P\cdot W)_{211213}=\mp \frac{4}{3}e^{2 X_3}(2b(X_{3})+b'(X_{3})),\nonumber\\
&(P\cdot W)_{121312}=-(P\cdot W)_{211312}=\mp \frac{2}{3}e^{2 X_3}(2b(X_{3})+b'(X_{3})),\nonumber\\
&(P\cdot W)_{121323}=\mp \frac{2}{3}+\frac{1}{2}e^{2 X_3}(2b(X_{3})+b'(X_{3}))^2,\nonumber\\
&(P\cdot W)_{122311}=-(P\cdot W)_{212311}=\pm 2e^{2 X_3}(2b(X_{3})+b'(X_{3})),\nonumber
\end{align}
\begin{align}
&(P\cdot W)_{122313}=-(P\cdot W)_{212313}=\pm \frac{2}{3},\nonumber\\
&(P\cdot W)_{122322}=-(P\cdot W)_{212322}=\pm 2e^{-2 X_3}(2b(X_{3})+b'(X_{3})),\nonumber\\
&(P\cdot W)_{122333}=-(P\cdot W)_{212333}=2(2b(X_{3})+b'(X_{3})),\nonumber\\
&(P\cdot W)_{131212}=-(P\cdot W)_{311212}=\mp \frac{2}{3}e^{2 X_3}(2b(X_{3})+b'(X_{3})),\nonumber\\
&(P\cdot W)_{131223}=-(P\cdot W)_{311223}=\mp \frac{2}{3},\nonumber\\
&(P\cdot W)_{131313}=\frac{1}{3}e^{4 X_3}(2b(X_{3})+b'(X_{3})),\nonumber\\
&(P\cdot W)_{132323}=-(P\cdot W)_{312323}=\frac{1}{3}(2b(X_{3})+b'(X_{3})),\nonumber\\
&(P\cdot W)_{211323}=\pm \frac{2}{3}+\frac{1}{2}e^{2 X_3}(2b(X_{3})+b'(X_{3}))^2,\nonumber\\
&(P\cdot W)_{221223}=(P\cdot W)_{222312}=\mp e^{-2 X_3}(2b(X_{3})+b'(X_{3})),\nonumber\\
&(P\cdot W)_{221313}=e^{2 X_3}(2b(X_{3})+b'(X_{3}))^2,\nonumber\\
&(P\cdot W)_{231211}=-(P\cdot W)_{321211}=\pm 2e^{2 X_3}(2b(X_{3})+b'(X_{3})),\nonumber\\
&(P\cdot W)_{231213}=-(P\cdot W)_{321213}=\pm \frac{2}{3},\nonumber\\
&(P\cdot W)_{231222}=-(P\cdot W)_{321222}=\pm 2e^{-2 X_3}(2b(X_{3})+b'(X_{3})),\nonumber\\
&(P\cdot W)_{231233}=-(P\cdot W)_{321233}=2(2b(X_{3})+b'(X_{3})),\nonumber\\
&(P\cdot W)_{231312}=(P\cdot W)_{321312}=\frac{1}{2}e^{2 X_3}(2b(X_{3})+b'(X_{3}))^2,\nonumber\\
&(P\cdot W)_{231323}=-\frac{2}{3}(2b(X_{3})+b'(X_{3})),\nonumber\\
&(P\cdot W)_{232313}=-(P\cdot W)_{322313}=\frac{2}{3}(2b(X_{3})+b'(X_{3})),\nonumber\\
&(P\cdot W)_{311313}=-\frac{7}{3}e^{4 X_3}(2b(X_{3})+b'(X_{3})),\nonumber\\
&(P\cdot W)_{321323}=-\frac{4}{3}(2b(X_{3})+b'(X_{3})),\nonumber\\
&(P\cdot W)_{331223}=(P\cdot W)_{332312}=-(2b(X_{3})+b'(X_{3})),\nonumber\\
&(P\cdot W)_{331313}=\pm 2e^{4 X_3}(2b(X_{3})+b'(X_{3}))^2;\nonumber\\
&(P\cdot P)_{111223}=-(P\cdot P)_{112123}=\mp e^{2 X_3}(2b(X_{3})+b'(X_{3})),\\
&(P\cdot P)_{111313}=-(P\cdot P)_{113113}=2e^{6 X_3}(2b(X_{3})+b'(X_{3}))^2,\nonumber\\
&(P\cdot P)_{112312}=-(P\cdot P)_{113212}=\mp e^{2 X_3}(2b(X_{3})+b'(X_{3})),\nonumber\\
&(P\cdot P)_{121123}=-(P\cdot P)_{211123}=\pm \frac{1}{2}e^{2 X_3}(2b(X_{3})+b'(X_{3})),\nonumber\\
&(P\cdot P)_{121213}=-(P\cdot P)_{211213}=\mp e^{2 X_3}(2b(X_{3})+b'(X_{3})),\nonumber
\end{align}
\begin{align}
&(P\cdot P)_{121312}=-(P\cdot P)_{211312}=\mp \frac{1}{2}e^{2 X_3}(2b(X_{3})+b'(X_{3})),\nonumber\\
&(P\cdot P)_{121323}=\mp 1+\frac{1}{2}e^{2 X_3}(2b(X_{3})+b'(X_{3}))^2,\nonumber\\
&(P\cdot P)_{122113}=-(P\cdot P)_{212113}=\pm e^{2 X_3}(2b(X_{3})+b'(X_{3})),\nonumber\\
&(P\cdot P)_{122223}=-(P\cdot P)_{212223}=\pm \frac{1}{2}e^{-2 X_3}(2b(X_{3})+b'(X_{3})),\nonumber\\
&(P\cdot P)_{122311}=-(P\cdot P)_{212311}=\pm 2e^{2 X_3}(2b(X_{3})+b'(X_{3})),\nonumber\\
&(P\cdot P)_{122313}=-(P\cdot P)_{212313}=\pm 1,\nonumber\\
&(P\cdot P)_{122322}=-(P\cdot P)_{212322}=\pm 2e^{-2 X_3}(2b(X_{3})+b'(X_{3})),\nonumber\\
&(P\cdot P)_{122333}=-(P\cdot P)_{212333}=2(2b(X_{3})+b'(X_{3})),\nonumber\\
&(P\cdot P)_{123123}=- \frac{3}{4}e^{2 X_3}(2b(X_{3})+b'(X_{3}))^2,\nonumber\\
&(P\cdot P)_{123211}=-(P\cdot P)_{213211}=\mp e^{2 X_3}(2b(X_{3})+b'(X_{3})),\nonumber\\
&(P\cdot P)_{123213}=-(P\cdot P)_{213213}=-\frac{1}{2}e^{2 X_3}(2b(X_{3})+b'(X_{3}))^2,\nonumber\\
&(P\cdot P)_{123222}=-(P\cdot P)_{213222}=\mp e^{-2 X_3}(2b(X_{3})+b'(X_{3})),\nonumber\\
&(P\cdot P)_{123233}=-(P\cdot P)_{213233}=-(2b(X_{3})+b'(X_{3})),\nonumber\\
&(P\cdot P)_{123312}=-(P\cdot P)_{213312}=-\frac{1}{4}e^{2 X_3}(2b(X_{3})+b'(X_{3}))^2,\nonumber\\
&(P\cdot P)_{123323}=-(P\cdot P)_{213323}=\frac{1}{2}(2b(X_{3})+b'(X_{3})),\nonumber\\
&(P\cdot P)_{131111}=-(P\cdot P)_{311111}=\mp e^{6 X_3}(2b(X_{3})+b'(X_{3})),\nonumber\\
&(P\cdot P)_{131113}=-(P\cdot P)_{311113}=\mp e^{4 X_3}- e^{6 X_3}(2b(X_{3})+b'(X_{3}))^2,\nonumber\\
&(P\cdot P)_{131122}=-(P\cdot P)_{311122}=\mp e^{2 X_3}(2b(X_{3})+b'(X_{3})),\nonumber\\
&(P\cdot P)_{131133}=-(P\cdot P)_{311133}=-e^{4 X_3}(2b(X_{3})+b'(X_{3})),\nonumber\\
&(P\cdot P)_{131212}=-(P\cdot P)_{311212}=\mp \frac{1}{2}e^{2 X_3}(2b(X_{3})+b'(X_{3})),\nonumber\\
&(P\cdot P)_{131223}=-(P\cdot P)_{311223}=\mp 1- \frac{1}{4}e^{2 X_3}(2b(X_{3})+b'(X_{3}))^2,\nonumber\\
&(P\cdot P)_{131313}=\frac{3}{2}e^{4 X_3}(2b(X_{3})+b'(X_{3})),\nonumber\\
&(P\cdot P)_{132112}=-(P\cdot P)_{312112}=\pm e^{2 X_3}(2b(X_{3})+b'(X_{3})),\nonumber\\
&(P\cdot P)_{132323}=-(P\cdot P)_{312323}=2b(X_{3})+b'(X_{3}),\nonumber\\
&(P\cdot P)_{133113}=\frac{3}{2}e^{4 X_3}(2b(X_{3})+b'(X_{3})),\nonumber\\
&(P\cdot P)_{133212}=-(P\cdot P)_{313212}=-\frac{1}{4}e^{2 X_3}(2b(X_{3})+b'(X_{3}))^2,\nonumber
\end{align}
\begin{align}
&(P\cdot P)_{133223}=-(P\cdot P)_{313223}=\frac{1}{2}(2b(X_{3})+b'(X_{3})),\nonumber\\
&(P\cdot P)_{133311}=-(P\cdot P)_{313311}=e^{4 X_3}(2b(X_{3})+b'(X_{3})),\nonumber\\
&(P\cdot P)_{133313}=-(P\cdot P)_{313313}=\mp e^{4 X_3}(2b(X_{3})+b'(X_{3}))^2,\nonumber\\
&(P\cdot P)_{133322}=-(P\cdot P)_{313322}=2b(X_{3})+b'(X_{3}),\nonumber\\
&(P\cdot P)_{133333}=-(P\cdot P)_{313333}=\pm e^{2 X_3}(2b(X_{3})+b'(X_{3})),\nonumber\\
&(P\cdot P)_{211323}=\pm 1+ \frac{1}{2}e^{2 X_3}(2b(X_{3})+b'(X_{3}))^2,\nonumber\\
&(P\cdot P)_{213123}=-\frac{1}{4}e^{2 X_3}(2b(X_{3})+b'(X_{3}))^2,\nonumber\\
&(P\cdot P)_{221223}=-(P\cdot P)_{222123}=\mp e^{-2 X_3}(2b(X_{3})+b'(X_{3})),\nonumber\\
&(P\cdot P)_{221313}=-(P\cdot P)_{223113}=e^{2 X_3}(2b(X_{3})+b'(X_{3}))^2,\nonumber\\
&(P\cdot P)_{222312}=-(P\cdot P)_{223212}=\mp e^{-2 X_3}(2b(X_{3})+b'(X_{3})),\nonumber\\
&(P\cdot P)_{231123}=-(P\cdot P)_{321123}=-\frac{1}{4}e^{2 X_3}(2b(X_{3})+b'(X_{3}))^2,\nonumber\\
&(P\cdot P)_{231211}=-(P\cdot P)_{321211}=\pm e^{2 X_3}(2b(X_{3})+b'(X_{3})),\nonumber\\
&(P\cdot P)_{231213}=-(P\cdot P)_{321213}=-\frac{1}{2}e^{2 X_3}(2b(X_{3})+b'(X_{3}))^2,\nonumber\\
&(P\cdot P)_{231222}=-(P\cdot P)_{321222}=\pm e^{-2 X_3}(2b(X_{3})+b'(X_{3})),\nonumber\\
&(P\cdot P)_{231233}=-(P\cdot P)_{321233}=2b(X_{3})+b'(X_{3}),\nonumber\\
&(P\cdot P)_{231312}=\frac{1}{4}e^{2 X_3}(2b(X_{3})+b'(X_{3}))^2,\nonumber\\
&(P\cdot P)_{231323}=-\frac{1}{2}(2b(X_{3})+b'(X_{3})),\nonumber\\
&(P\cdot P)_{232111}=-(P\cdot P)_{322111}=\mp 2e^{2 X_3}(2b(X_{3})+b'(X_{3})),\nonumber\\
&(P\cdot P)_{232113}=-(P\cdot P)_{322113}=\mp 1,\nonumber\\
&(P\cdot P)_{232122}=-(P\cdot P)_{322122}=\mp 2e^{-2 X_3}(2b(X_{3})+b'(X_{3})),\nonumber\\
&(P\cdot P)_{232133}=-(P\cdot P)_{322133}=-2(2b(X_{3})+b'(X_{3})),\nonumber\\
&(P\cdot P)_{232212}=-(P\cdot P)_{322212}=\pm \frac{1}{2}e^{-2 X_3}(2b(X_{3})+b'(X_{3})),\nonumber\\
&(P\cdot P)_{232223}=-(P\cdot P)_{322223}=\mp e^{-4 X_3},\nonumber\\
&(P\cdot P)_{232313}=-(P\cdot P)_{322313}=2b(X_{3})+b'(X_{3}),\nonumber\\
&(P\cdot P)_{233112}=(P\cdot P)_{323112}=-\frac{1}{2}e^{2 X_3}(2b(X_{3})+b'(X_{3}))^2,\nonumber\\
&(P\cdot P)_{233123}=(P\cdot P)_{323123}=2b(X_{3})+b'(X_{3}),\nonumber\\
&(P\cdot P)_{233312}=-(P\cdot P)_{323312}=2b(X_{3})+b'(X_{3}),\nonumber
\end{align}
\begin{align}
&(P\cdot P)_{311313}=-\frac{7}{2}e^{4 X_3}(2b(X_{3})+b'(X_{3})),\nonumber\\
&(P\cdot P)_{313113}=\frac{1}{2}e^{4 X_3}(2b(X_{3})+b'(X_{3})),\nonumber\\
&(P\cdot P)_{321312}=\frac{3}{4}e^{2 X_3}(2b(X_{3})+b'(X_{3}))^2,\nonumber\\
&(P\cdot P)_{321323}=-\frac{3}{2}(2b(X_{3})+b'(X_{3})),\nonumber\\
&(P\cdot P)_{331223}=-(P\cdot P)_{332123}=-(2b(X_{3})+b'(X_{3})),\nonumber\\
&(P\cdot P)_{331313}=-(P\cdot P)_{333113}=\pm 2e^{4 X_3}(2b(X_{3})+b'(X_{3}))^2,\nonumber\\
&(P\cdot P)_{332312}=-(P\cdot P)_{333212}=-(2b(X_{3})+b'(X_{3})).\nonumber
\end{align}
By calculations, we have
\begin{align}
&Q(Ric, R)_{111313}=-4e^{6 X_3}(2b(X_{3})+b'(X_{3}))^2,\\
&Q(Ric, R)_{121213}=-Q(Ric, R)_{211213}=\pm 2e^{2 X_3}(2b(X_{3})+b'(X_{3})),\nonumber\\
&Q(Ric, R)_{121312}=-Q(Ric, R)_{211312}=\pm e^{2 X_3}(2b(X_{3})+b'(X_{3})),\nonumber\\
&Q(Ric, R)_{121323}=\mp 2-e^{2 X_3}(2b(X_{3})+b'(X_{3}))^2,\nonumber\\
&Q(Ric, R)_{122313}=-Q(Ric, R)_{212313}=\pm 2,\nonumber\\
&Q(Ric, R)_{131212}=-Q(Ric, R)_{311212}=\pm e^{2 X_3}(2b(X_{3})+b'(X_{3})),\nonumber\\
&Q(Ric, R)_{131223}=-Q(Ric, R)_{311223}=\mp 2,\nonumber\\
&Q(Ric, R)_{131313}=-2e^{4 X_3}(2b(X_{3})+b'(X_{3})),\nonumber\\
&Q(Ric, R)_{132323}=-Q(Ric, R)_{312323}=-(2b(X_{3})+b'(X_{3})),\nonumber\\
&Q(Ric, R)_{211323}=\pm 2-e^{2 X_3}(2b(X_{3})+b'(X_{3}))^2,\nonumber\\
&Q(Ric, R)_{221313}=-2e^{2 X_3}(2b(X_{3})+b'(X_{3}))^2,\nonumber\\
&Q(Ric, R)_{231213}=-Q(Ric, R)_{321213}=\pm 2,\nonumber\\
&Q(Ric, R)_{231312}=-Q(Ric, R)_{321312}=-e^{2 X_3}(2b(X_{3})+b'(X_{3}))^2,\nonumber\\
&Q(Ric, R)_{231323}=2b(X_{3})+b'(X_{3}),\nonumber\\
&Q(Ric, R)_{232313}=-Q(Ric, R)_{322313}=-2(2b(X_{3})+b'(X_{3})),\nonumber\\
&Q(Ric, R)_{311313}=6e^{4 X_3}(2b(X_{3})+b'(X_{3})),\nonumber\\
&Q(Ric, R)_{321323}=3(2b(X_{3})+b'(X_{3})),\nonumber\\
&Q(Ric, R)_{331313}=\mp 4e^{4 X_3}(2b(X_{3})+b'(X_{3}))^2;\nonumber\\
&Q(Ric, C)_{111313}=-4e^{6 X_3}(2b(X_{3})+b'(X_{3}))^2,\\
&Q(Ric, C)_{121323}=-2e^{2 X_3}(2b(X_{3})+b'(X_{3}))^2,\nonumber\\
&Q(Ric, C)_{122313}=-Q(Ric, C)_{212313}=-2e^{2 X_3}(2b(X_{3})+b'(X_{3}))^2,\nonumber
\end{align}
\begin{align}
&Q(Ric, C)_{131223}=-Q(Ric, C)_{311223}=e^{2 X_3}(2b(X_{3})+b'(X_{3}))^2,\nonumber\\
&Q(Ric, C)_{131313}=Q(Ric, C)_{311313}=2e^{4 X_3}(2b(X_{3})+b'(X_{3})),\nonumber\\
&Q(Ric, C)_{132312}=-Q(Ric, C)_{312312}=-e^{2 X_3}(2b(X_{3})+b'(X_{3}))^2,\nonumber\\
&Q(Ric, C)_{132323}=-Q(Ric, C)_{312323}=2(2b(X_{3})+b'(X_{3})),\nonumber\\
&Q(Ric, C)_{221313}=-2e^{2 X_3}(2b(X_{3})+b'(X_{3}))^2,\nonumber\\
&Q(Ric, C)_{231213}=-Q(Ric, C)_{321213}=2e^{2 X_3}(2b(X_{3})+b'(X_{3}))^2,\nonumber\\
&Q(Ric, C)_{321312}=-2e^{2 X_3}(2b(X_{3})+b'(X_{3}))^2,\nonumber\\
&Q(Ric, C)_{321323}=4(2b(X_{3})+b'(X_{3})),\nonumber\\
&Q(Ric, C)_{331313}=\mp 4e^{4 X_3}(2b(X_{3})+b'(X_{3}))^2;\nonumber\\
&Q(Ric, K)_{111313}=-4e^{6 X_3}(2b(X_{3})+b'(X_{3}))^2,\\
&Q(Ric, K)_{121213}=-Q(Ric, K)_{211213}=\pm 2e^{2 X_3}(2b(X_{3})+b'(X_{3})),\nonumber\\
&Q(Ric, K)_{121312}=-Q(Ric, K)_{211312}=\pm e^{2 X_3}(2b(X_{3})+b'(X_{3})),\nonumber\\
&Q(Ric, K)_{121323}=\mp 2-2e^{2 X_3}(2b(X_{3})+b'(X_{3}))^2,\nonumber\\
&Q(Ric, K)_{122313}=-Q(Ric, K)_{212313}=\pm 2-2e^{2 X_3}(2b(X_{3})+b'(X_{3}))^2,\nonumber\\
&Q(Ric, K)_{131212}=-Q(Ric, K)_{311212}=\pm e^{2 X_3}(2b(X_{3})+b'(X_{3})),\nonumber\\
&Q(Ric, K)_{131223}=-Q(Ric, K)_{311223}=\mp 2+e^{2 X_3}(2b(X_{3})+b'(X_{3}))^2,\nonumber\\
&Q(Ric, K)_{131313}=6e^{4 X_3}(2b(X_{3})+b'(X_{3})),\nonumber\\
&Q(Ric, K)_{132312}=-Q(Ric, K)_{312312}=-e^{2 X_3}(2b(X_{3})+b'(X_{3}))^2,\nonumber\\
&Q(Ric, K)_{132323}=-Q(Ric, K)_{312323}=3(2b(X_{3})+b'(X_{3})),\nonumber\\
&Q(Ric, K)_{211323}=\pm 2,\nonumber\\
&Q(Ric, K)_{221313}=-2e^{2 X_3}(2b(X_{3})+b'(X_{3}))^2,\nonumber\\
&Q(Ric, K)_{231213}=-Q(Ric, K)_{321213}=\pm 2+2e^{2 X_3}(2b(X_{3})+b'(X_{3}))^2,\nonumber\\
&Q(Ric, K)_{231323}=2b(X_{3})+b'(X_{3}),\nonumber\\
&Q(Ric, K)_{232313}=-Q(Ric, K)_{322313}=2(2b(X_{3})+b'(X_{3})),\nonumber\\
&Q(Ric, K)_{311313}=-2e^{4 X_3}(2b(X_{3})+b'(X_{3})),\nonumber\\
&Q(Ric, K)_{321312}=-2e^{2 X_3}(2b(X_{3})+b'(X_{3}))^2,\nonumber\\
&Q(Ric, K)_{321323}=3(2b(X_{3})+b'(X_{3})),\nonumber\\
&Q(Ric, K)_{331313}=\mp 4e^{4 X_3}(2b(X_{3})+b'(X_{3}))^2;\nonumber\\
&Q(Ric, W)_{111313}=-4e^{6 X_3}(2b(X_{3})+b'(X_{3}))^2,\\
&Q(Ric, W)_{121213}=-Q(Ric, W)_{211213}=\pm \frac{8}{3}e^{2 X_3}(2b(X_{3})+b'(X_{3})),\nonumber\\
&Q(Ric, W)_{121312}=-Q(Ric, W)_{211312}=\pm \frac{4}{3}e^{2 X_3}(2b(X_{3})+b'(X_{3})),\nonumber\\
&Q(Ric, W)_{121323}=\mp \frac{8}{3}-e^{2 X_3}(2b(X_{3})+b'(X_{3}))^2,\nonumber
\end{align}
\begin{align}
&Q(Ric, W)_{122313}=-Q(Ric, W)_{212313}=\pm \frac{8}{3},\nonumber\\
&Q(Ric, W)_{131212}=-Q(Ric, W)_{311212}=\pm \frac{4}{3}e^{2 X_3}(2b(X_{3})+b'(X_{3})),\nonumber\\
&Q(Ric, W)_{131223}=-Q(Ric, W)_{311223}=\mp \frac{8}{3},\nonumber\\
&Q(Ric, W)_{131313}=-\frac{2}{3}e^{4 X_3}(2b(X_{3})+b'(X_{3})),\nonumber\\
&Q(Ric, W)_{132323}=-Q(Ric, W)_{312323}=-\frac{2}{3}(2b(X_{3})+b'(X_{3})),\nonumber\\
&Q(Ric, W)_{211323}=\pm \frac{8}{3}-e^{2 X_3}(2b(X_{3})+b'(X_{3}))^2,\nonumber\\
&Q(Ric, W)_{221313}=-2e^{2 X_3}(2b(X_{3})+b'(X_{3}))^2,\nonumber\\
&Q(Ric, W)_{231213}=-Q(Ric, W)_{321213}=\pm \frac{8}{3},\nonumber\\
&Q(Ric, W)_{231312}=Q(Ric, W)_{321312}=-e^{2 X_3}(2b(X_{3})+b'(X_{3}))^2,\nonumber\\
&Q(Ric, W)_{231323}=\frac{4}{3}(2b(X_{3})+b'(X_{3})),\nonumber\\
&Q(Ric, W)_{232313}=-Q(Ric, W)_{322313}=-\frac{4}{3}(2b(X_{3})+b'(X_{3})),\nonumber\\
&Q(Ric, W)_{311313}=\frac{14}{3}e^{4 X_3}(2b(X_{3})+b'(X_{3})),\nonumber\\
&Q(Ric, W)_{321323}=\frac{8}{3}(2b(X_{3})+b'(X_{3})),\nonumber\\
&Q(Ric, W)_{331313}=\mp 4e^{4 X_3}(2b(X_{3})+b'(X_{3}))^2;\nonumber\\
&Q(Ric, P)_{111313}=-Q(Ric, P)_{113113}=-4e^{6 X_3}(2b(X_{3})+b'(X_{3}))^2,\\
&Q(Ric, P)_{121213}=-Q(Ric, P)_{211213}=\pm 2e^{2 X_3}(2b(X_{3})+b'(X_{3})),\nonumber\\
&Q(Ric, P)_{121312}=-Q(Ric, P)_{211312}=\pm e^{2 X_3}(2b(X_{3})+b'(X_{3})),\nonumber\\
&Q(Ric, P)_{121323}=\mp 2- e^{2 X_3}(2b(X_{3})+b'(X_{3}))^2,\nonumber\\
&Q(Ric, P)_{122113}=-Q(Ric, P)_{212113}=\mp 2e^{2 X_3}(2b(X_{3})+b'(X_{3})),\nonumber\\
&Q(Ric, P)_{122313}=-Q(Ric, P)_{212313}=\pm 2,\nonumber\\
&Q(Ric, P)_{123112}=-Q(Ric, P)_{213112}=\mp e^{2 X_3}(2b(X_{3})+b'(X_{3})),\nonumber\\
&Q(Ric, P)_{123123}=\pm 2+ \frac{3}{2}e^{2 X_3}(2b(X_{3})+b'(X_{3}))^2,\nonumber\\
&Q(Ric, P)_{123213}=-Q(Ric, P)_{213213}=\mp 2+e^{2 X_3}(2b(X_{3})+b'(X_{3}))^2,\nonumber\\
&Q(Ric, P)_{123312}=-Q(Ric, P)_{213312}=\frac{1}{2}e^{2 X_3}(2b(X_{3})+b'(X_{3}))^2,\nonumber\\
&Q(Ric, P)_{123323}=-Q(Ric, P)_{213323}=-(2b(X_{3})+b'(X_{3})),\nonumber\\
&Q(Ric, P)_{131113}=-Q(Ric, P)_{311113}=2e^{6 X_3}(2b(X_{3})+b'(X_{3}))^2,\nonumber
\end{align}
\begin{align}
&Q(Ric, P)_{131212}=-Q(Ric, P)_{311212}=\pm e^{2 X_3}(2b(X_{3})+b'(X_{3})),\nonumber\\
&Q(Ric, P)_{131223}=-Q(Ric, P)_{311223}=\mp 2+\frac{1}{2}e^{2 X_3}(2b(X_{3})+b'(X_{3}))^2,\nonumber\\
&Q(Ric, P)_{131313}=-3e^{4 X_3}(2b(X_{3})+b'(X_{3})),\nonumber\\
&Q(Ric, P)_{132112}=-Q(Ric, P)_{312112}=\mp e^{2 X_3}(2b(X_{3})+b'(X_{3})),\nonumber\\
&Q(Ric, P)_{132123}=-Q(Ric, P)_{312123}=\pm 2,\nonumber\\
&Q(Ric, P)_{132323}=-Q(Ric, P)_{312323}=-(2b(X_{3})+b'(X_{3})),\nonumber\\
&Q(Ric, P)_{133113}=-3e^{4 X_3}(2b(X_{3})+b'(X_{3})),\nonumber\\
&Q(Ric, P)_{133212}=-Q(Ric, P)_{313212}=\frac{1}{2}e^{2 X_3}(2b(X_{3})+b'(X_{3}))^2,\nonumber\\
&Q(Ric, P)_{133223}=-Q(Ric, P)_{313223}=-(2b(X_{3})+b'(X_{3})),\nonumber\\
&Q(Ric, P)_{133313}=-Q(Ric, P)_{313313}=2e^{2 X_3}\pm 2e^{4 X_3}(2b(X_{3})+b'(X_{3}))^2,\nonumber\\
&Q(Ric, P)_{211323}=\pm 2- e^{2 X_3}(2b(X_{3})+b'(X_{3}))^2,\nonumber\\
&Q(Ric, P)_{213123}=\mp 2+ \frac{1}{2}e^{2 X_3}(2b(X_{3})+b'(X_{3}))^2,\nonumber\\
&Q(Ric, P)_{221313}=-Q(Ric, P)_{223113}=-2e^{2 X_3}(2b(X_{3})+b'(X_{3}))^2,\nonumber\\
&Q(Ric, P)_{231123}=-Q(Ric, P)_{321123}=\frac{1}{2}e^{2 X_3}(2b(X_{3})+b'(X_{3}))^2,\nonumber\\
&Q(Ric, P)_{231213}=-Q(Ric, P)_{321213}=\pm 2+e^{2 X_3}(2b(X_{3})+b'(X_{3}))^2,\nonumber\\
&Q(Ric, P)_{231312}=-\frac{1}{2}e^{2 X_3}(2b(X_{3})+b'(X_{3}))^2,\nonumber\\
&Q(Ric, P)_{232113}=-Q(Ric, P)_{322113}=\mp 2,\nonumber\\
&Q(Ric, P)_{232313}=-Q(Ric, P)_{322313}=-2(2b(X_{3})+b'(X_{3})),\nonumber\\
&Q(Ric, P)_{233112}=e^{2 X_3}(2b(X_{3})+b'(X_{3}))^2,\nonumber\\
&Q(Ric, P)_{233123}=-2(2b(X_{3})+b'(X_{3})),\nonumber\\
&Q(Ric, P)_{233312}=-Q(Ric, P)_{323312}=-(2b(X_{3})+b'(X_{3})),\nonumber\\
&Q(Ric, P)_{233323}=-Q(Ric, P)_{323323}=2e^{-2 X_3},\nonumber\\
&Q(Ric, P)_{311313}=7e^{4 X_3}(2b(X_{3})+b'(X_{3})),\nonumber\\
&Q(Ric, P)_{313113}=-e^{4 X_3}(2b(X_{3})+b'(X_{3})),\nonumber\\
&Q(Ric, P)_{321312}=-\frac{3}{2}e^{2 X_3}(2b(X_{3})+b'(X_{3}))^2,\nonumber\\
&Q(Ric, P)_{321323}=4(2b(X_{3})+b'(X_{3})),\nonumber\\
&Q(Ric, P)_{323112}=e^{2 X_3}(2b(X_{3})+b'(X_{3}))^2,\nonumber\\
&Q(Ric, P)_{323123}=-2(2b(X_{3})+b'(X_{3})),\nonumber\\
&Q(Ric, P)_{331313}=-Q(Ric, P)_{333113}=\mp 4e^{4 X_3}(2b(X_{3})+b'(X_{3}))^2.\nonumber
\end{align}
It follows immediately that
\begin{align}
&Q(g, R)_{111223}=-Q(g, R)_{112312}=-e^{2 X_3}(2b(X_{3})+b'(X_{3})),\\
&Q(g, R)_{121323}=Q(g, R)_{131223}=-2,\nonumber\\
&Q(g, R)_{122313}=Q(g, R)_{231213}=2,\nonumber\\
&Q(g, R)_{221223}=-Q(g, R)_{222312}=-e^{-2 X_3}(2b(X_{3})+b'(X_{3})),\nonumber\\
&Q(g, R)_{331223}=-Q(g, R)_{332312}=\mp (2b(X_{3})+b'(X_{3}));\nonumber\\
&Q(g, C)_{111223}=-Q(g, C)_{112312}=-e^{2 X_3}(2b(X_{3})+b'(X_{3})),\\
&Q(g, C)_{121312}=-Q(g, C)_{131212}=-e^{2 X_3}(2b(X_{3})+b'(X_{3})),\nonumber\\
&Q(g, C)_{132323}=-Q(g, C)_{231323}=\pm (2b(X_{3})+b'(X_{3})),\nonumber\\
&Q(g, C)_{221223}=-Q(g, C)_{222312}=-e^{-2 X_3}(2b(X_{3})+b'(X_{3})),\nonumber\\
&Q(g, C)_{331223}=-Q(g, C)_{332312}=\mp (2b(X_{3})+b'(X_{3}));\nonumber\\
&Q(g, K)_{111223}=-Q(g, K)_{112312}=-e^{2 X_3}(2b(X_{3})+b'(X_{3})),\\
&Q(g, K)_{121312}=-Q(g, K)_{131212}=-e^{2 X_3}(2b(X_{3})+b'(X_{3})),\nonumber\\
&Q(g, K)_{132323}=-Q(g, K)_{231323}=\pm (2b(X_{3})+b'(X_{3})),\nonumber\\
&Q(g, K)_{221223}=-Q(g, K)_{222312}=-e^{-2 X_3}(2b(X_{3})+b'(X_{3})),\nonumber\\
&Q(g, K)_{331223}=-Q(g, K)_{332312}=\mp (2b(X_{3})+b'(X_{3}));\nonumber\\
&Q(g, W)_{111223}=-Q(g, W)_{112312}=-e^{2 X_3}(2b(X_{3})+b'(X_{3})),\\
&Q(g, W)_{121323}=Q(g, W)_{131223}=-2,\nonumber\\
&Q(g, W)_{122313}=Q(g, W)_{231213}=2,\nonumber\\
&Q(g, W)_{221223}=-Q(g, W)_{222312}=-e^{-2 X_3}(2b(X_{3})+b'(X_{3})),\nonumber\\
&Q(g, W)_{331223}=-Q(g, W)_{332312}=\mp (2b(X_{3})+b'(X_{3}));\nonumber\\
&Q(g, P)_{111223}=-Q(g, P)_{112123}=-e^{2 X_3}(2b(X_{3})+b'(X_{3})),\\
&Q(g, P)_{112312}=-Q(g, P)_{113212}=e^{2 X_3}(2b(X_{3})+b'(X_{3})),\nonumber\\
&Q(g, P)_{121123}=\frac{1}{2}e^{2 X_3}(2b(X_{3})+b'(X_{3})),\nonumber\\
&Q(g, P)_{121323}=-2,~Q(g, P)_{123123}=1,\nonumber\\
&Q(g, P)_{122223}=\frac{1}{2}e^{-2 X_3}(2b(X_{3})+b'(X_{3})),\nonumber\\
&Q(g, P)_{122313}=2,~Q(g, P)_{123213}=-1,\nonumber\\
&Q(g, P)_{123112}=\frac{1}{2}e^{2 X_3}(2b(X_{3})+b'(X_{3})),\nonumber
\end{align}
\begin{align}
&Q(g, P)_{131113}=-e^{4 X_3},~Q(g, P)_{133313}=\pm e^{2 X_3},\nonumber\\
&Q(g, P)_{131223}=-2,~Q(g, P)_{132123}=1,\nonumber\\
&Q(g, P)_{132112}=-\frac{1}{2}e^{2 X_3}(2b(X_{3})+b'(X_{3})),\nonumber\\
&Q(g, P)_{132323}=\pm \frac{1}{2}(2b(X_{3})+b'(X_{3})),\nonumber\\
&Q(g, P)_{221223}=-Q(g, P)_{222312}=-e^{-2 X_3}(2b(X_{3})+b'(X_{3})),\nonumber\\
&Q(g, P)_{231213}=1,~Q(g, P)_{232113}=-2,\nonumber\\
&Q(g, P)_{231323}=\mp \frac{1}{2}(2b(X_{3})+b'(X_{3})),\nonumber\\
&Q(g, P)_{232212}=-\frac{1}{2}e^{-2 X_3}(2b(X_{3})+b'(X_{3})),\nonumber\\
&Q(g, P)_{232223}=-e^{-4 X_3},~Q(g, P)_{233323}=\pm e^{-2 X_3},\nonumber\\
&Q(g, P)_{233312}=\mp \frac{1}{2}(2b(X_{3})+b'(X_{3})),\nonumber\\
&Q(g, P)_{331223}=-Q(g, P)_{332312}=\mp (2b(X_{3})+b'(X_{3})).\nonumber
\end{align}

On account of the above result, we have
\begin{thm}
The metric $(1.1)$ fulfills the following geometric structures about the connection $\widehat{\nabla}$:

(I)$R\cdot K=0$ from this it is the conharmonic curvature semisymmetric type manifold due to the Riemann curvature.

(II)The condition $Ric^{3}\mp 2 Ric^{2}=0$ holds, hence it is a Einstein manifold of level $3.$

(III)Ricci tensor is neither of Codazzi type nor cyclic parallel.

(IV)The Riemann curvature, the Weyl conformal curvature, the conharmonic curvature, the concircular curvature and the projective curvature are not recurrent.
\end{thm}
\vskip 1 true cm
\section{The geometric structures in $Sol_{3}$ with the Levi-Civita connection}
According to the definition of the semi-symmetric non-metric connection, when $a(X_{3})=0$ or $b(X_{3})=0$ we can get the Levi-Civita connection $\nabla=\widehat{\nabla}.$   Meanwhile, the Christofel coefficient of the Levi-Civita connection $\Gamma^{\alpha}_{kj}$ are given below
\begin{align}
&\Gamma^{1}_{13}=\Gamma^{1}_{31}=1,~\Gamma^{2}_{23}=\Gamma^{2}_{32}=-1,~\Gamma^{3}_{11}=\mp e^{2 X_3},~\Gamma^{3}_{22}=\pm e^{-2 X_3}.
\end{align}
It is clear that
\begin{align}
&R_{1212}=\pm 1,~R_{1313}=-e^{2 X_3},~R_{2323}=-e^{-2 X_3};~Ric_{33}=2;~\kappa=\pm 2;\\
&\nabla_{2} R_{1213}=\pm 2,~\nabla_{1} R_{1223}=\pm 2;~\nabla_{1} Ric_{13}=\pm 2e^{2 X_3},~\nabla_{2} Ric_{23}=\mp 2e^{-2 X_3}.
\end{align}
Likewise,
\begin{align}
&K_{1212}=\pm 1,~K_{1313}=e^{2 X_3},~K_{2323}=e^{-2 X_3};\\
&W_{1212}=\pm \frac{4}{3},~W_{1313}=-\frac{2}{3}e^{2 X_3},~W_{2323}=-\frac{2}{3}e^{-2 X_3};\\
&P_{1212}=-P_{1221}=\pm 1,~P_{1313}=-e^{2 X_3},~P_{2323}=-e^{-2 X_3}.
\end{align}
Then, we obtain
\begin{align}
&\nabla_{2} W_{1213}=\nabla_{2} W_{1312}=\pm 2,~\nabla_{1} W_{1223}=\nabla_{1} W_{2312}=\pm 2;\\
&\nabla_{2} P_{1213}=\nabla_{2} P_{1312}=\pm 2,~\nabla_{1} P_{1223}=-\nabla_{1} P_{2321}=\pm 2,\\
&\nabla_{2} P_{1231}=\nabla_{2} P_{1321}=\mp 1,~\nabla_{1} P_{1232}=-\nabla_{1} P_{2312}=\mp 1,\nonumber\\
&\nabla_{1} P_{1311}=\mp e^{4 X_3},~\nabla_{1} P_{1333}=e^{2 X_3},~\nabla_{2} P_{2322}=\pm e^{-4 X_3},~\nabla_{2} P_{2333}=-e^{-2 X_3}.\nonumber
\end{align}

Similar to Theorem 3.1, we can get the following theorem
\begin{thm}
The metric $(1.1)$ fulfills the following geometric structures about the connection $\nabla$:

(I)(1)$C\cdot R=0$ from this it is the Riemann curvature semisymmetric type manifold due to the Weyl conformal curvature;\\
(2)$K\cdot R=-Q(Ric, R)$ hence it is Ricci generalized the Riemann curvature pseudosymmetric type manifold due to the conharmonic curvature;\\
(3)$W\cdot R=\frac{2}{3}Q(Ric, R)$ thus it is Ricci generalized the Riemann curvature pseudosymmetric type manifold due to the concircular curvature;\\
(4)Ricci generalized projectively pseudosymmetric because $P\cdot R=\frac{1}{2}Q(Ric, R)$;\\
(5)$K\cdot R=\mp Q(g, R)$ so it is the Riemann curvature pseudosymmetric type manifold due to the conharmonic curvature;\\
(6)$W\cdot R=\pm \frac{2}{3}Q(g, R)$ from this it is the Riemann curvature pseudosymmetric type manifold due to the concircular curvature;\\
(7)$P\cdot R=\pm \frac{1}{2}Q(g, R)$ i.e. it is the Riemann curvature pseudosymmetric type manifold due to the projective curvature;\\
(8)$R\cdot C=0$ then it is the Weyl conformal curvature semisymmetric type manifold due to the Riemann curvature;\\
(9)$C\cdot C=0$ for this reason it is the Weyl conformal curvature semisymmetric type manifold due to the Weyl conformal curvature;\\
(10)$K\cdot C=0$ in this way it is the Weyl conformal curvature semisymmetric type manifold due to the conharmonic curvature;\\
(11)$W\cdot C=0$ hence it is the Weyl conformal curvature semisymmetric type manifold due to the concircular curvature;\\
(12)$P\cdot C=0$ so that it is the Weyl conformal curvature semisymmetric type manifold due to the projective curvature;\\
(13)$R\cdot K=0$ i.e. it is the conharmonic curvature semisymmetric type manifold due to the Riemann curvature;\\
(14)$C\cdot K=0$ for this reason it is the conharmonic curvature semisymmetric type manifold due to the Weyl conformal curvature;\\
(15)$K\cdot K=0$ so that it is the conharmonic curvature semisymmetric type manifold due to the conharmonic curvature;\\
(16)$W\cdot K=0$ in this way it is the conharmonic curvature semisymmetric type manifold due to the concircular curvature;\\
(17)$P\cdot K=-\frac{1}{2}Q(Ric, K)$ for this reason it is Ricci generalized the conharmonic curvature pseudosymmetric type manifold due to the projective curvature;\\
(18)$R\cdot W=\frac{3}{4}Q(Ric, W)$ hence it is Ricci generalized the concircular curvature pseudosymmetric type manifold due to the Riemann curvature;\\
(19)$K\cdot W=-\frac{3}{4}Q(Ric, W)$ from this it is Ricci generalized the concircular curvature pseudosymmetric type manifold due to the conharmonic curvature;\\
(20)$W\cdot W=\frac{1}{2}Q(Ric, W)$ therefore it is Ricci generalized the concircular curvature pseudosymmetric type manifold due to the concircular curvature;\\
(21)$P\cdot W=\frac{1}{4}Q(Ric, W)$ so it is Ricci generalized the concircular curvature pseudosymmetric type manifold due to the projective curvature;\\
(22)$R\cdot W=\pm Q(g, W)$ then it is the concircular curvature pseudosymmetric type manifold due to the Riemann curvature;\\
(23)Concircular curvature pseudosymmetric type manifold due to the conharmonic curvature i.e. $K\cdot W=\mp Q(g, W)$;\\
(24)$W\cdot W=\pm \frac{2}{3}Q(g, W)$ in this way it is the concircular curvature pseudosymmetric type manifold due to the concircular curvature;\\
(25)$P\cdot W=\pm \frac{1}{3}Q(g, W)$ thus it is the concircular curvature pseudosymmetric type manifold due to the projective curvature;\\
(26)$R\cdot P=\pm Q(g, P)$ therefore it is the projective curvature pseudosymmetric type manifold due to the Riemann curvature;\\
(27)$K\cdot P=\mp Q(g, P)$ for this reason it is the projective curvature pseudosymmetric type manifold due to the conharmonic curvature;\\
(28)$W\cdot P=\pm \frac{2}{3}Q(g, P)$  hence it is the projective curvature pseudosymmetric type manifold due to the concircular curvature.

(II)$Ric=\alpha(\eta\bigotimes\eta)$ for $\alpha=2$ and $\eta=\{0, 0, 1\},$ thus it is a Ricci simple manifold and a $2$-quasi-Einstein manifold as $rank(Ric-\alpha g)=2$ when $g_{33}=1,$ or a $3$-quasi-Einstein manifold as $rank(Ric-\alpha g)=3$ when $g_{33}=-1.$

(III)The manifold is the generalized Roter type manifold since $R=g\wedge Ric\mp \frac{1}{2} g\wedge g.$

(IV)The condition $Ric^{2}\mp 2 Ric=0$ holds, hence it is a Einstein manifold of level $2.$

(V)Ricci tensor is neither of Codazzi type nor cyclic parallel.

(VI)Ricci tensor is Riemann compatible, Weyl conformal compatible, conharmonic compatible and concircular compatible.
\end{thm}

\vskip 1 true cm

%-----------------------------------------------------------------------------
%-----------------------------------------------------------------------------

\bigskip
\bigskip

\noindent {\footnotesize {\it S. Liu} \\
{School of Mathematics and Statistics, Northeast Normal University, Changchun 130024, China}\\
{Email: liusy719@nenu.edu.cn}

\end{document}